\documentclass[12pt]{article}
\usepackage[latin1]{inputenc}
\usepackage{bbm}
\usepackage{latexsym}
\usepackage{amsfonts}
\usepackage{amsmath,amssymb,amscd}
\usepackage{stmaryrd}

\DeclareFontFamily{U}{txsyc}{}
\DeclareFontShape{U}{txsyc}{m}{n}{
   <-> txsyc
}{}
\DeclareFontShape{U}{txsyc}{bx}{n}{
   <-> txbsyc
}{}
\DeclareFontShape{U}{txsyc}{l}{n}{<->ssub * txsyc/m/n}{}
\DeclareFontShape{U}{txsyc}{b}{n}{<->ssub * txsyc/bx/n}{}

\DeclareSymbolFont{symbolsC}{U}{txsyc}{m}{n}
\SetSymbolFont{symbolsC}{bold}{U}{txsyc}{bx}{n}
\DeclareFontSubstitution{U}{txsyc}{m}{n}

\DeclareMathSymbol{\df}{\mathrel}{symbolsC}{"42}
\DeclareMathSymbol{\fd}{\mathrel}{symbolsC}{"43}
\DeclareMathSymbol{\lJoin}{\mathrel}{symbolsC}{"58}
\DeclareMathSymbol{\rJoin}{\mathrel}{symbolsC}{"59}

\newcommand{\cA}{{\cal A}}

\newcommand{\cD}{{\cal D}}

\newcommand{\cF}{{\cal F}}

\newcommand{\cK}{{\cal K}}
\newcommand{\cL}{{\cal L}}
\newcommand{\cN}{{\cal N}}
\newcommand{\cM}{{\cal M}}
\newcommand{\cO}{{\cal O}}
\newcommand{\cP}{{\cal P}}

\newcommand{\cT}{{\cal T}}

\newcommand{\CC}{\mathbb{C}}

\newcommand{\EE}{\mathbb{E}}

\newcommand{\LL}{\mathbb{L}}
\newcommand{\NN}{\mathbb{N}}
\newcommand{\PP}{\mathbb{P}}

\newcommand{\RR}{\mathbb{R}}

\newcommand{\di}{\displaystyle}
\newcommand{\iy}{\infty}
\newcommand{\lt}{\left}
\newcommand{\me}{\medskip}

\newcommand{\ri}{\rightarrow}
\newcommand{\rt}{\right}
\newcommand{\sm}{\smallskip}
\newcommand{\tr}{\triangle}
\newcommand{\wi}{\widetilde}

\newcommand{\card}{\mathrm{card}}

\newcommand{\diam}{\mathrm{diam}}

\newcommand{\ex}{\exists\ }
\newcommand{\fo}{\forall\ }

\newcommand{\lve}{\lt\vert}
\newcommand{\lVe}{\lt\Vert}

\newcommand{\rve}{\rt\vert}
\newcommand{\rVe}{\rt\Vert}

\newcommand{\st}{\,:\,}

\newcommand{\un}{1{\hskip -2.5 pt}\hbox{I}}

\newcommand{\bq}{\begin{eqnarray*}}
\newcommand{\bqn}[1]{\begin{eqnarray}\label{#1}}
\newcommand{\eq}{\end{eqnarray*}}
\newcommand{\eqn}{\end{eqnarray}}

\newcommand{\wwtbp}{\par\hfill $\blacksquare$\par\me\noindent}

\newcommand{\thistitlepagestyle}{}

\newcommand{\lin}{\llbracket}
\newcommand{\rin}{\rrbracket}

\newcommand{\ttsim}{\raise.17ex\hbox{$\scriptstyle\mathtt{\sim}$}}

\newcommand{\thh}[1]{${#1}^{\mathrm{th}}$}

\newtheorem{pro}{Proposition} 

\newtheorem{lem}[pro]{Lemma}
\newtheorem{theo}[pro]{Theorem}
\renewcommand{\thepro}{\arabic{pro}}

\newenvironment{rem}
{\par\me\refstepcounter{pro}\noindent{\bf Remark \thepro\ }}
{\par\hfill $\square$\par\me\noindent}

\newcommand{\proof}{\par\me\noindent\textbf{Proof}\par\sm\noindent}

\title{On barycentric subdivision, with simulations}
\author{Persi Diaconis${}^\dagger$ and Laurent Miclo${}^\ddagger$}
\date{\box1
\vskip-2mm
\box2
\vskip3mm
\box3}

\begin{document}

\setbox1=\vbox{
\large
\begin{center}
${}^\dagger$Departments of Statistics and Mathematics\\
Stanford University, USA\\
\end{center}
} 
\setbox2=\vbox{
\large
\begin{center}
${}^\ddagger$Institut de Mathématiques de Toulouse, UMR 5219\\
Universit\'e de Toulouse and CNRS, France
\end{center}
} 
\setbox3=\vbox{
\Large
\begin{center}
We dedicate this paper to the memory of David Blackwell.
\end{center}
}
\setbox4=\vtop{
\hbox{${}^\dagger$diaconis@math.stanford.edu \\}
\vskip1mm
\hbox{Department of Statistics\\}
\hbox{Sequoia Hall\\}
\hbox{390 Serra Mall\\}
\hbox{Stanford University\\}
\hbox{Stanford, CA 94305-4065\\}
}
\setbox5=\vtop{
\hbox{${}^\ddagger$ miclo@math.univ-toulouse.fr\\}
\vskip1mm
\hbox{Institut de Mathématiques de Toulouse\\}
\hbox{Université Paul Sabatier\\}
\hbox{118, route de Narbonne\\} 
\hbox{31062 Toulouse Cedex 9, France\\}
}
\setbox6=\hbox{
\box4\hskip20mm
\box5
}

\maketitle
\thistitlepagestyle
\abstract{Consider the barycentric subdivision which cuts a given triangle along its medians
to produce six new triangles. Uniformly choosing one of them and iterating this procedure gives rise
to a Markov chain. We show that almost surely, the triangles forming this chain become flatter and flatter
in the sense that their isoperimetric values goes to infinity with time. Nevertheless, if the triangles are renormalized
through a similitude to have their longest edge equal to $[0,1]\subset\CC$ (with 0 also adjacent to the shortest edge), their aspect does not converge and we identify the
limit set of the opposite vertex with the segment [0,1/2]. In addition we  prove that the largest angle converges to $\pi$ in probability. Our approach is probabilistic and these results are deduced from the investigation of a limit iterated random function Markov chain living on the segment [0,1/2].
The stationary distribution of this limit chain is particularly important in our study. 
In an appendix we present related numerical simulations (not included in
the version submitted for publication).
}
\vfill\null
{\small
\textbf{Keywords: }barycentric subdivision, triangle-valued Markov chain, isoperimetric functional, flat triangles,
iterated random functions, invariant probability measure, Wasserstein distance.
\par
\vskip.3cm
\textbf{MSC2000:} first: 60J05,
secondary: 
60D05, 60F15, 26A18, 28C10, 37A30.
}\par

\newpage

\section{Introduction}

Let $\tr$ be a given triangle in the plane (to avoid triviality the vertices will always be assumed not to be all the same). The three medians of $\tr$ intersect at the barycenter, this cuts it into six small triangles,
say $\tr_1,\ \tr_2,\ \tr_3,\ \tr_4,\ \tr_5,\ \tr_6$.
Next, each $\tr_i$, for $i\in\lin 1, 6\rin$ (which  denotes the set $\{1,2, ..., 6\}$), can itself be subdivided in the same way into 
six  triangles, $(\tr_{i,j})_{j\in\lin 1, 6\rin}$. Iterating this barycentric subdivision procedure,
we get $6^n$ triangles $(\tr_{I})_{I\in\lin 1, 6\rin^n}$ at stage $n\in\NN$. It is well-known numerically 
(we learned it from Blackwell \cite{Blackwell_barycentric}, see also the survey by Butler and Graham 
\cite{Butler_triangle})
and  has been recently proved (cf.\ Diaconis and McMullen \cite{Diaconis_barycentric} and Hough \cite{MR2516262}) that as the barycentric subdivision goes on, most of the triangles  become flat. The original motivation for this kind of result was to show that the barycentric subdivision is not a good procedure to construct nice triangularizations of surfaces. For more information on other kinds of triangle
subdivisions, we refer to a recent manuscript of Butler and Graham 
\cite{Butler_triangle}.
The goal of this paper is to propose a new probabilistic approach to this phenomenon.
\par\me
First, we adopt a Markovian point of view: Let $\tr(0)\df\tr$ and throw a fair die to choose $\tr(1)$ among the six
triangles $\tr_i$, $i\in\lin 1, 6\rin$. Continuing in the same way, we get a Markov chain $(\tr(n))_{n\in\NN}$: 
if the \thh{n} first triangles have been constructed, the next one is obtained
by choosing uniformly (and independently from what was done before) one of the six triangles of the barycentric subdivision
of the last obtained triangle. Of course, at any time $n\in\NN^*$ ($\NN^*$ stands for $\NN\setminus\{0\}$), the law of $\tr(n)$ is the uniform distribution on the set of triangles
$\{\tr_{I}\st I\in\lin 1, 6\rin^n\}$. So to deduce generic properties under this distribution it is sufficient to study the chain
$(\tr(n))_{n\in\NN}$.\\
In order to describe our results more analytically, let us renormalize the triangles. For any non-trivial triangle
$\tr$ on the plane, there is a similitude of the plane transforming $\tr$ into a triangle whose vertices
are $(0,0)$, $(0,1)$ and $(x,y)\in[0,1/2]\times [0,\sqrt{3}/2]$, such that the longest (respectively the shortest) edge of $\tr$
is sent to $[(0,0),  (0,1)]$ (resp.\ $[(0,0),(x,y)]$). The point $(x,y)$ is uniquely determined and characterizes the aspect of $\tr$ (as long as orientation is not considered, otherwise we would have to consider positive similitudes and $x$ would have to belong to $[0,1]$). 
Any time we are interested in quantities which are invariant by similitude, we will identify triangles with
their characterizing points. In particular, this identification will endow the set of triangles with the topology
(not separating triangles with the same aspect) inherited from the usual topology of the plane. This convention will 
implicitly be enforced  throughout this paper.
The triangle $\tr$ will be said to be flat if $y=0$. So up to similitude the set of flat triangles can be identified with
$[0,1/2]$.
For $n\in\NN$, let $(X_n,Y_n)$  be the characterizing point of $\tr(n)$. 
The first result  justifies  the assertion that as the barycentric subdivision goes on, the triangles become flat.
\begin{theo}\label{th1}
Almost surely (a.s.) the stochastic sequence $(Y_n)_{n\in\NN}$ converges to zero exponentially fast: there exists a constant $\chi>0$ such that a.s.:
\bq \limsup_{n\ri\iy}\frac1{n}\ln(Y_n)&\leq&-\chi\eq
\end{theo}
We will show that we can take $\chi= 0.035$ 
(but this is not the best constant, indeed we will present a numerical experiment suggesting that the above bound should hold with $\chi\approx 0.07$), nevertheless the previous result remains asymptotical. 
But contrary to Blackwell \cite{Blackwell_barycentric} (see also the remark at the end of section 6),
we have not been able to deduce a more quantitative bound in probability on $Y_n$ for any given $n\in\NN$.
\par\sm
In particular, we recover the convergence in probability toward the set of flat triangles which was previously proved by Diaconis and McMullen \cite{Diaconis_barycentric} (using Furstenberg theorem on 
products of random matrices in $SL_2(\RR)$)
and Hough \cite{MR2516262} who used
dynamical systems arguments (via an identification with a random walk on $SL_2(\RR)$). 
\par\sm
There is a stronger notion of convergence to flatness that asks for the triangles to have an angle which is almost equal to $\pi$.
With the preceding notation, for $n\in\NN$, let $A_n$ be the angle between $[(0,0),(X_n,Y_n)]$ and $[(X_n,Y_n),(0,1)]$, this is the largest
angle of $\tr(n)$.
\begin{theo}\label{th2}
The sequence $(\tr(n))_{n\in\NN}$ becomes
strongly  flat in probability:
\bq
\fo \epsilon>0,\qquad \lim_{n\ri\iy}\PP[A_n<\pi-\epsilon]&=&0\eq
\end{theo}\par
Of course this result implies that $(Y_n)_{n\in\NN}$ converges to zero in probability. Note that the converse is not true in general: there are isosceles triangles that become flatter and flatter, but their maximum
angle converges to $\pi/2$.
Indeed, Theorem \ref{th2} is more difficult to obtain than Theorem~\ref{th1}
because $(X_n)_{n\in\NN}$ does not converge as the next result shows. Define the limit set of this sequence as the intersection over $p\in\NN$
of the closures of the sets $\{X_n\st n\geq p\}$.
\begin{theo}\label{th3}
Almost surely, the limit set of $(X_n)_{n\in\NN}$ is $[0,1/2]$. 
\end{theo}
It follows from Theorem \ref{th1}
that a.s.\ the limit set of a trajectory of the triangle Markov chain $(\tr(n))_{n\in\NN}$
is the whole set of flat triangles.\par\sm
A crucial tool behind these results is a limiting flat Markov chain $Z$.
 Strictly speaking the stochastic chain $(X_n)_{n\in\NN}$ is not Markovian, but eventually its evolution becomes almost Markovian. Indeed, we note that the above barycentric subdivision procedure can formally also be applied to flat triangles
and their set is stable by this operation. This means that if $Y_0=0$, then for any $n\in\NN$, $Y_n=0$ a.s. 
In this particular situation $(X_n)_{n\in\NN}$ is Markovian. Let $M$ be its transition kernel, from $[0,1/2]$ to itself. In what follows,
$Z\df(Z_n)_{n\in\NN}$ will always designate a Markov chain on $[0,1/2]$ whose transition kernel is $M$.
An important part of this paper will be devoted to the investigation of the Markov chain $Z$ since it is the key to the above asymptotic behaviour.
We will see that $Z$ is ergodic in the sense that it admits an attracting (and thus unique) invariant measure $\mu$ on $[0,1/2]$. We will also show that $\mu$ is continuous and that its support is $[0,1/2]$ 
(but we don't know if $\mu$ is absolutely continuous).
\par\me
The plan of the paper is the following: the next section 
contains some global preliminaries, in particular we will show, by
studying the evolution of a convenient variant of the isoperimetric value, that the triangle Markov chain
 returns as  close as we want to the set of flat triangles infinitely often. This is a first step in the direction of Theorem \ref{th1}.
In section 3, 
we begin our investigation of the limiting Markov chain $Z$, to obtain some information valid in a neighborhood of the set of flat triangles.
Then in section 4 we put together the previous global and local results to prove Theorem \ref{th1}.
Ergodicity and the attracting invariant measure $\mu$ of the Markov chain $Z$ are studied in section 5, using results of 
Dubins and Freedman \cite{MR0193668},
Barnsley and Elton \cite{MR932532} 
and Diaconis and Freedman \cite{MR1669737}
on 
iterated random functions. This will lead to the proofs of Theorem \ref{th2}  and Theorem \ref{th3} in section~6.
In an appendix we present some numerical experiments which helped us in some technical proofs.
 They will also illustrate the obtained results, in particular to evaluate the constant $\chi$ of Theorem~\ref{th1} and to see the profile of $\mu$. 

\section{A weak result on attraction to flatness}

The purpose of this section is to give some preliminary information and bounds on the triangle Markov chain obtained by barycentric subdivisions.
By themselves, these results are not sufficient to conclude the a.s.\ convergence toward the set of flat triangles, but at least
they give a heuristic hint for this behaviour: 
a quantity comparable to
the isoperimetry value of the triangle has a tendency to increase after barycentric subdivision
and so to diverge to infinity with time, in the mean.
\par\me
To measure the separation between a given triangle $\tr$
 and the set of flat triangles, we use the quantity $J(\tr)$, which is the sum of the squares of the lengths of the edges divided by the area (this is well-defined in $(0,+\iy]$, since the vertices are assumed not to be all the same).
 We have $J(\tr)=+\iy$ if and only if $\tr$ is flat. Furthermore the functional $J$ is invariant under similitude,
 so it depends only on the characteristic point $(x,y)$ of $\tr$ and we have
 \bq
 J(\tr)\ =\ {2}\frac{1+{x^2+y^2}+{(1-x)^2+y^2}}{{y}}\ =\ 4\frac{x^2+y^2-x+1}{y}\eq
 in particular we get
 \bqn{yI}
 3y\ \leq\ (J(\tr))^{-1}\ \leq\ 8y\eqn 
 so that the convergence of $y$ to zero is equivalent to the divergence of $J(\tr)$ to $+\iy$.
 Note that $J(\tr)$ is comparable with the isoperimetric value $I(\tr)$ of $\tr$, defined as the square of the perimeter of $\tr$ divided by its area of $\tr$:
\bqn{IJ}
\frac13I(\tr)\ \leq\ J(\tr)\ \leq\ I(\tr)\eqn
 With the notation of the introduction, write 
 for $n\in\NN$, $J_n\df J(\tr(n))$. Our first goal is to show
 \begin{pro}\label{div}
 Almost surely, we have $\limsup_{n\ri\iy}J_n=+\iy$.
 \end{pro}
 The proof will be based on elementary considerations of one step of the barycentric subdivision.
 Consider $\tr$ a triangle in the normalized form given in the introduction.
 For simplicity, we denote $A$, $B$ and $C$ the vertices $(0,0)$, $(x,y)$ and $(1,0)$ of $\tr$.
 Let also $D$, $E$, $F$ and $G$ be respectively the middle points of $[A,B]$, $[B,C]$
 and $[A,C]$ and the barycenter of $\tr$. 
 We index the small triangles obtained by the barycentric subdivision as
 \bqn{triangles}
 \begin{array}{ccc}
 \tr_1\ \df\ \{A,D,G\},\qquad&
 \tr_2\ \df\ \{D,B,G\},\qquad&
 \tr_3\ \df\ \{B,E,G\}\\
 \tr_4\ \df\ \{E,C,G\},\qquad&
 \tr_5\ \df\ \{C,F,G\},\qquad&
 \tr_6\ \df\ \{F,A,G\}
 \end{array}\eqn\par
It is well-known that all the triangles $\tr_i$, for $i\in\lin 1, 6\rin$, have the same area
(one straightforward way to see it is to note that this property is invariant by
affine transformations and to consider the equilateral case).\par
Next  define, with $\lve\cdot\rve$ denoting the length,
\bq
L_1\ \df\ \lve[A,B]\rve,\qquad
L_2\ \df\ \lve[B,C]\rve,\qquad
L_3\ \df\ \lve[C,A]\rve\\
l_1\ \df\ \lve[D,C]\rve,\qquad
l_2\ \df\ \lve[E,A]\rve,\qquad
l_3\ \df\ \lve[F,B]\rve\eq
\par
An immediate computation gives that
\bqn{lll}
l_1^2\ =\ \frac{x^2}4+\frac{y^2}4-x+1,\qquad
l_2^2\ =\ \frac{x^2}4+\frac{y^2}4+\frac{x}2+\frac14,\qquad
l_3^2\ =\ {x^2}+{y^2}-x+\frac14\eqn
so we get
\bq
\frac{l_1^2+l_2^2+l_3^2}{L_1^2+L_2^2+L_3^2}
&=&\frac34\eq
These ingredients imply the following probabilistic statement:
\begin{lem}\label{subm}
For any $n\in\NN$, we have  
\bq
\EE[J_{n+1}\vert \cT_n]&=&\frac43J_n\eq
where the lhs is a conditional expectation with respect to $\cT_n$, the $\sigma$-algebra generated by $\tr(n)$, $\tr(n-1)$, ..., $\tr(0)$.
\end{lem}
\proof
By the Markov property, the above bound is equivalent to the fact that for any $n\in\NN$,
\bq
\EE[J_{n+1}\vert \tr(n)]&=&\frac43J_n\eq
Since the Markov chain $(\tr(n))_{n\in\NN}$ is time-homogeneous, it is sufficient
to deal with the case $n=0$.
We come back to the notation introduced above. Because the small triangles have the same area
and the barycenter cuts the median segments into a ratio (1/3,2/3),
we get that
\bq
\EE[J(\tr(1))\vert \tr(0)=\tr]&=&\frac16\sum_{i\in\lin1,6\rin}J(\tr_i)\\
&=&\frac16\lt(\frac{L_1^2}2+\frac{L_2^2}2+\frac{L_3^2}2+\frac{10}{9}(l_1^2+l_2^2+l_3^2)\rt)
\frac6{\cA(\tr)}\\
&=& \frac43J(\tr)\eq 
where $\cA(\tr)$ is the area of $\tr$.
\wwtbp\par
In general the previous submartingale information is not enough to deduce a.s.\ convergence.
Taking expectations, we get that for any $n\in\NN$,
$
\EE[J_{n+1}]\geq(4/3)\EE[J_n]$, 
thus
$\EE[J_n]\geq (4/3)^nJ(\tr)$,
so we can just deduce $\LL^1$-divergence of $J_n$ for large $n\in\NN$, but this is not a very useful result.\par\sm
To prove Proposition \ref{div}, note that the numbers $J_n$, $n\in\NN$, are uniformly bounded below by a positive constant.
Indeed the isoperimetric functional can be defined for relatively general bounded subsets of the plane 
and its minimal value is obtained for discs (cf.\ for instance Osserman \cite{MR0500557}). Thus, via (\ref{IJ}), we get
\bq
\fo n\in\NN, \qquad J_n\ \geq\ \frac43\pi\ >\ 4\eq
\par
But from Lemma \ref{subm}, we see that
\bq
\fo n\in\NN, \qquad \PP[J_{n+1}\geq(4/3)J_n\vert \cT_n]&\geq&\frac16\eq
and consequently
\bqn{plonge}
\fo n,m\in\NN, \qquad \PP[J_{n+m}\geq(4/3)^m4\vert \cT_n]&\geq&\frac1{6^m}\eqn
Let $R>1$ be an arbitrary large number and consider $m\in\NN^*$  such that $(4/3)^m4\geq R$.
The $\{0,1\}$-valued sequence $(\un_{J_{m(n+1)}\geq R})_{n\in\NN}$ stochastically dominates a sequence of independent Bernoulli variables of parameter $1/6^m$.
It follows that a.s.\ we have
\bq
\limsup_{n\ri\iy}J_n&\geq& R\eq
and since $R$ can be chosen arbitrarily large, Proposition \ref{div} is proven.\wwtbp
\par\me
To finish this section, we will prove another simple preliminary result.
\begin{lem}\label{ab}
There exists two constants $0<a<b<+\iy$ such that
\bq
\fo n\in\NN, \qquad aJ_n\ \leq\ J_{n+1}\ \leq\ bJ_n\eq
\end{lem}
\proof
Again it is sufficient to consider the first barycentric subdivision and to prove
that we can find two constants $0<a<b<+\iy$ such that with the above notation,
\bq
\fo i\in\lin 1,6\rin,\qquad aJ(\tr)\ \leq\ J(\tr_i)\ \leq\ bJ(\tr)\eq
Such inequalities are obvious for flat triangles, so  assume that $\tr$ is not flat.
Since the areas are easy to compare, we just need to consider the 
diameters (whose squares are comparable, within the range [1,3], 
with the sums of the squares of the lengths of the edges),
denoted by $d$.
We have clearly $d(\tr)=1$ and $d(\tr_i)\leq1$ for $i\in\lin 1,6\rin$.
The reverse bound $d(\tr_i)\geq 1/4$, for $i\in\lin 1,6\rin$,
is a consequence of the equalities $\lve [A,F]\rve=\lve[F,C]\rve=1/2$,
$\lve [B,E]\rve=\lve [E,C]\rve\geq1/2$
and $\lve [D,G]\rve=\lve [G,C]\rve/2\geq1/4$.
\wwtbp

\section{Near the limit flat Markov chain}

Our goal here is two-fold. First we show that the kernel of the triangle Markov chain
converges nicely to the kernel of the flat triangle Markov chain as the triangle
becomes flat. Second we study the evolution of a perimeter related functional
for the flat triangle Markov chain, to get a bound on the evolution of the isoperimetric functional
for the triangle Markov chain, valid at least in a neighborhood of the set of flat triangles.\par\me
Let $Q$ be the transition kernel of the Markov chain $(X_n,Y_n)_{n\in\NN}$ considered in the introduction.
For any $(x,y)\in\cD$, the set of characterizing points of triangles, we can write
\bq
Q((x,y),\,\cdot\,)&=&\frac16\sum_{i\in\lin 1,6\rin}\delta_{(x_i,y_i)}\eq 
where $\delta$ stands for the Dirac mass and where for any $i\in\lin 1,6\rin$, $(x_i,y_i)$ is the characterizing point of the triangle $\tr_i$ described in (\ref{triangles}). Of course, the $x_i$ and $y_i$, for $i\in\lin 1,6\rin$, have to be seen as functions of $(x,y)$.
For $i\in\lin 1,6\rin$, let us define
\bqn{zx}
\fo x\in [0,1/2],\qquad z_i(x)&\df& x_i(x,0)\eqn 
The transition kernel $M$ on $[0,1/2]$ of the flat triangle Markov chain alluded to in the introduction can expressed as
\bqn{Mz}\fo x\in[0,1/2],\qquad
M(x,\,\cdot\,)&=&\frac16\sum_{i\in\lin 1,6\rin}\delta_{z_i(x)}\eqn 
The next result gives bounds on the discrepancy between $Q$ and $M$ as the triangles become flat.
\begin{lem}\label{QM}
There exists a constant $K>0$ such that
\bq
\fo i\in\lin 1,6\rin,\,\fo (x,y)\in\cD,\qquad 
\max(
\lve x_i(x,y)-z_i(x)\rve,
\lve y_i(x,y)\rve)&\leq& Ky
\eq
\end{lem}
\proof
We first check that for any fixed $i\in\lin 1,6\rin$,
the mapping
\bqn{xiyi}
\cD_2\ni(x,y^2)&\mapsto&(x_i(x,y),y_i^2(x,y))\in\cD_2\eqn
is (uniformly) Lipschitz, where $\cD_2$ is the image
of $\cD$ through $\cD\ni (x,y)\mapsto (x,y^2)\in\cD_2$.
\par
Indeed, denote by $0\leq L_{i,1}\leq L_{i,2}\leq L_{i,3}$ the ordered
lengths of the triangle $\tr_i$.
We have seen in the previous section that
\bqn{yi}
y_i^2&=&\frac{L_{i,1}^2+ L_{i,2}^2+ L_{i,3}^2}{2L_{i,3}^2}-x_i^2+x_i-1\eqn
Let $h_i$ be the height of $\tr_i$ orthogonal to the edge of length $L_{i,3}$, we have
$
L_{i,1}^2=h_i^2+(x_iL_{i,3})^2$ and $
L_{i,2}^2=h_i^2+((1-x_i)L_{i,3})^2$. 
It follows that 
\bqn{xi}
x_i&=&\frac{L_{i,3}^2-L_{i,2}^2+L_{i,1}^2}{2L_{i,3}^2}\eqn
Finally, notice that (\ref{lll}) implies that the mappings
$\cD_2\ni(x,y^2)\mapsto L_{i,j}^2$, for $j \in\lin 1,3\rin$, are uniformly Lipschitz.
Furthermore, as seen in the proof of Lemma \ref{ab},
the mapping $\cD_2\ni(x,y^2)\mapsto L_{i,3}^2$ is bounded below by 1/16,
so (\ref{xi}) and (\ref{yi}) imply that the mapping described in (\ref{xiyi}) is
uniformly Lipschitz.\par
The bounds given in Lemma \ref{QM} are an easy consequence 
of this Lipschitz property and of the boundedness of $\cD$.\wwtbp\par
The second goal of this section is to study the sign of quantities like
$
\EE[\ln(I_{n+1}/I_n)\vert \tr(n)=\tr]$, at least when $\tr$ is close to a flat triangle. Here we define $I_n$ as the isoperimetric value of $\tr(n)$. 
By the Markov property, this amounts to evaluating the sign of
$1/6\sum_{i\in\lin 1,6\rin}\ln(I(\tr_i)/I(\tr))$.
Of course the previous ratios are not rigorously defined if the triangle $\tr$ is flat, 
Nevertheless, let $(x,y)$ be the characterizing point of $\tr$. When $y$ goes to zero $0_+$, 
$I(\tr_i)/I(\tr)=\sqrt{6}\cP(\tr_i)/\cP(\tr)$ converges to $G(i,x)$, which is just the same ratio for the flat triangle $\tr$ whose
characterizing point is $(x,0)$. We have, for any $x\in[0,1/2]$ (see the computations of section 5 for more details)
\bqn{G}\begin{array}{rclcrcl}
G(1,x)&=&\sqrt{\frac23}(1+x),&\qquad&
G(2,x)&=&\sqrt{\frac16}(2-x)\cr
G(3,x)&=&\sqrt{\frac32}(1-x),&\qquad&
G(4,x)&=&\sqrt{\frac23}(2-x)\cr
G(5,x)&=&\sqrt{\frac23}(2-x),&\qquad&
G(6,x)&=&\sqrt{\frac32}\end{array}\eqn
From the previous considerations, we easily get that this convergence is uniform over $x$, in the sense that
for any $i\in\lin 1,6\rin$,
\bq
\lim_{y\ri 0_+}\sup_{(x,y)\in\cD}\lve \frac{I(\tr_i)}{I(\tr)}-G(i,x)\rve&=&0\eq
So to prove that $
\EE[\ln(I_{n+1}/I_n)\vert \tr(n)=\tr]>0$ for nearly flat triangles $\tr$, it is sufficient to show
that the mapping $[0,1/2]\ni x\mapsto \sum_{i\in\lin 1,6\rin}\ln(G(i,x))$ only takes positive values.
Unfortunately, this is not true, since it takes negative values in a neighborhood of $1/2$ (see section 7).
To get around this problem, we iterate the barycentric subdivision one more step.
\begin{pro}\label{neigh}
There exist a constant $\gamma>0$ and a neighborhood $\cN$ of the set of the flat triangles, such that
\bq
\fo n\in\NN,\,\fo \tr\in \cN,\qquad \EE[\ln(I_{n+2}/I_n)\vert \tr(n)=\tr]&\geq&\gamma\eq
(for flat triangles $\tr$, the ratio is defined as a limit as above, or equivalently, as a ratio of perimeters, before renormalisation, up to the factor 6).
\end{pro}
\proof
Coming back to the notation at the beginning of the introduction, we want to find 
 $\cN$ and $\gamma$ as above and satisfying
 \bqn{gamma}
 \fo \tr\in \cN,\qquad \frac1{36}\sum_{i,j\in\lin 1,6\rin}\ln\lt(\frac{6\cP(\tr_{i,j})}{\cP(\tr)}\rt)&\geq&\gamma\eqn
 Let $(x,y)$ be the characterizing point of $\tr$. As $y$ goes to $0_+$, the left hand side
 converges (uniformly over $x$) to
 \bqn{F1}
 F(x)&\df& \frac1{36}\sum_{i,j\in\lin 1,6\rin}\ln(G(j,z_i(x))G(i,x))\eqn
 where the $z_i(x)$, for $i\in\lin 1,6\rin$, were defined in (\ref{zx}).
 More explicitly, we will compute in section 5 (see Lemma \ref{homz}) that on each of the segments
 $[0,1/5]$, $[1/5,2/7]$
 and $[2/7,1/2]$, the $z_i$, for $i\in\lin 1,6\rin$, are homographical mappings.
 So it seems more convenient to 
 consider
the piecewise rational fraction 
\bqn{R}
R(x)&\df&\exp(36 F(x))\\
\nonumber&=&\prod_{i,j\in\lin 1,6\rin}G(i,z_j(x))G(j,x)\eqn
After computations (see section 7), it appears that this is indeed 
a piecewise polynomial function.
By numerically studying the zeroes of $R-1$ of the
three underlying polynomial functions, 
we show that $F$ does not vanish on $[0,1/2]$.
So by continuity, we get that $\gamma\df\min_{[0,1/2]}F/2>0$. 
Then using the above uniform convergence, we can find a neighborhood $\cN$ of the set of flat triangles
so that (\ref{gamma}) is fulfilled.\wwtbp
We will see more precisely in section 7 that $F$ is decreasing, so we can take $\gamma=F(1/2)/2\approx 0.035$.

\section{Almost sure convergence to flatness}

We are now in position to prove Theorem \ref{th1}. The principle behind the proof is that
there is a neighborhood $\cN'$ of the set of flat triangles such that if the triangle Markov chain 
is inside $\cN'$, then it has a positive probability to always stay in this neighborhood and then
to converge exponentially fast to the set of flat triangles. This event will eventually occur, since
triangle Markov chains always return to $\cN'$.\par\me
In order to see that the triangle Markov chain has a positive chance to remain trapped in a neighborhood 
of the set of flat triangles, we will use a general martingale argument.
To do so, we introduce some notation.
On some underlying probability space, let $(\cF_n)_{n\in\NN}$ be a filtration, namely a non-decreasing sequence of 
$\sigma$-algebras. Let $\gamma>0$ and $A>0$ be two given constants.
We assume that for any $R$ large enough, say $R\geq R_0>0$, we are given a chain $(V^{(R)}_n)_{n\in\NN}$ and a martingale 
$(N^{(R)}_n)_{n\in\NN}$, adapted to the filtration $(\cF_n)_{n\in\NN}$, satisfying $V_0^{(R)}=R$,  
$N_0^{(R)}=0$ and such that for any time $n\in\NN$,
\bqn{NR}
\lve N^{(R)}_{n+1}-N^{(R)}_n\rve&\leq& A\\
\label{VR}
V^{(R)}_{n+1}-V^{(R)}_n&\geq&\gamma+ N^{(R)}_{n+1}-N^{(R)}_n\eqn
The next result shows that if $R$ is large enough, with high probability
$V^{(R)}$ will never go below $R/2$. This is classical, but without a precise reference at hand, we recall the underlying arguments.
\begin{lem}\label{convexp}
We have
\bq
\PP[\ex n\in\NN\st V_n^{(R)}<R/2]&\leq& \exp(-\gamma R/(2A^2))\frac{1}{1-\exp(-\gamma^2/(2A))}\eq
and furthermore,  a.s.,
\bq
\liminf_{n\ri\iy}\frac{V_n^{(R)}}{n}&\geq&\gamma\eq
\end{lem}
\proof
The first estimate is an immediate consequence of the Hoeffding-Azuma inequality,
which, applied to the  bounded difference martingale $(-N^{(R)}_n)_{n\in\NN}$ starting from 0,
asserts that for any $t\in\RR_+$,
\bq
\fo n\in\NN^*,\qquad\PP[-N^{(R)}_n>t]&\leq& \exp(-t^2/(2nA^2))\eq
In particular, since for any $n\in\NN$, we have
\bqn{VN}
V_n^{(R)}&\geq& R+n\gamma+N^{(R)}_n\eqn
we get
\bq
\PP[V_n^{(R)}<R/2]&\leq&\PP[-N^{(R)}_n>R/2+n\gamma]\\
&\leq& \exp\lt(-\frac{R}{4nA^2}-\frac{R\gamma}{2A^2}-\frac{n\gamma^2}{2A^2}\rt)\\
&\leq&\exp\lt(-\frac{R\gamma}{2A^2}\rt)\exp\lt(-\frac{n\gamma^2}{2A^2}\rt)\eq
and the first announced bound follows by summation over $n\in\NN^*$.
\par
The second bound is also due to the fact that the increments of the martingale $N^{(R)}$ are bounded, which implies the validity of the iterated logarithm law (see
for instance Stout \cite{MR0455094}): a.s.,
\bq
\limsup_{n\ri\iy}\frac{\lve N_n^{(R)}\rve}{\sqrt{n\ln(\ln(n))}}&\leq& A\eq
thus (\ref{VN}) enables us to conclude.
\wwtbp
Lemma \ref{convexp} will be applied with $V^{(R)}$ the logarithm of isoperimetric values, or rather with
a sequence of the kind $(\ln(I_{2n}))_{n\in\NN}$.
\\
More precisely, consider the neighborhhood $\cN$ obtained in Proposition \ref{neigh}. 
There exists a small constant $\epsilon>0$ such that
$\cN$ contains $\{(x,y)\in\cD\st 0\leq y<\epsilon\}$ and so taking into account (\ref{yI}), there exists
$R_1>1$ such that $\{\tr\st \ln(I(\tr))> R_1\}\subset \cN$ (again we are slightly abusing notation here,
identifying triangles with the characterizing points of their normalized forms, this should not lead to confusion).
Let $T$ be a finite stopping time for the triangle Markov chain $(\tr(n))_{n\in\NN}$. 
Assume that $R\df \ln(I(\tr(T)))$ satisfies $R\geq 2R_1$.
Define a stopping time $\tau$ for the shifted chain $(\tr(T+2n))_{n\in\NN}$ by 
\bq
\tau&\df&\inf\{n\in\NN\st \ln(I(\tr(T+2n)))\leq R_1\}\eq
which is infinite if the set on the r.h.s.\ is empty.
Let $\gamma>0$ be the constant appearing in Proposition~\ref{neigh}.
We construct a stochastic chain $V^{(R)}$ in the following way:
\bq
\fo n\in\NN,\qquad V_n^{(R)}&\df&\lt\{\begin{array}{ll}
\ln(I(\tr(T+2n)))&\hbox{, if $n\leq \tau$}\\
\ln(I(\tr(T+2\tau)))+\gamma(n-\tau)&\hbox{, otherwise}\end{array}\rt.\eq 
Let us check that the assumptions for Lemma \ref{convexp} are satisfied. 
Following the traditional Doob-Meyer semi-martingale decomposition (see for instance Dellacherie and Meyer \cite{MR566768}), we define 
\bq
\fo n\in\NN,\qquad N_n^{(R)}&\df&\sum_{m\in\lin 1,n\rin}V_m^{(R)}-\EE[V_m^{(R)}\vert \cF_{m-1}]\eq
where for any $n\in\NN$, $\cF_n$ is the $\sigma$-algebra generated by the trajectory-valued variable
$(\tr(m\wedge(T+n)))_{m\in\NN}$. Using classical stopping time notation, this is the $\sigma$-algebra
$\cT_{T+n}$, where the filtration $(\cT_m)_{m\in\NN}$ was introduced in Lemma \ref{subm}.
After conditioning on $\cF_0$ and taking advantage of the strong Markov property, 
we can apply Lemma \ref{ab} to see that (\ref{NR}) is satisfied with $A=(b/a)^2$ (we even have
$N^{(R)}_{n+1}-N^{(R)}_{n}=0$ for $n\geq \tau$).
Furthermore, we have for any $n\in\NN$,
\bq
V^{(R)}_{n+1}-V^{(R)}_n&=&\EE[V^{(R)}_{n+1}\vert\cF_n]-V^{(R)}_{n}+V^{(R)}_{n+1}-\EE[V^{(R)}_{n+1}\vert\cF_n]\\
&=&\EE[V^{(R)}_{n+1}-V^{(R)}_{n}\vert\cF_n]+N^{(R)}_{n+1}-N^{(R)}_n\\
&=&\EE[\ln(I_{T+2(n+1)}/I_{T+2n})\vert \tr(T+2n)]\un_{n\leq\tau}+\gamma\un_{n>\tau}+N^{(R)}_{n+1}-N^{(R)}_n
\\&\geq& \gamma+N^{(R)}_{n+1}-N^{(R)}_n\eq
where the last inequality comes from Proposition \ref{neigh}. Then Lemma \ref{convexp} implies that
\begin{pro}\label{stop}
Let $\cN'\df\{\tr\st \ln(I(\tr))>R_1\}$. There exists a large enough constant $R_2\geq2R_1$ such that
for any finite stopping time $T$  for the triangle Markov chain $(\tr(n))_{n\in\NN}$ satisfying
$\ln(I(\tr(T)))\geq R_2$, we have
\bq
\PP[\ex n\in\NN\st \tr(T+n)\not\in\cN'\vert \cT_{T}]&<&1/2\eq
Furthermore on the event $\{\fo n\in\NN\st \tr(T+n)\in\cN'\}$, we have a.s.
\bq
\liminf_{n\ri\iy}\frac{\ln(I_n)}{n}&\geq& \gamma/2\eq
\end{pro}
Indeed, Lemma \ref{convexp} shows that we can find  $R_2\geq2R_1$ such that
\bq
\PP[\tau<\iy\vert \cT_{T}]\ =\ \PP[\ex n\in\NN\st \tr(T+2n)\not\in\cN'\vert \cT_{T}]&<&1/2\eq
On the event $\{\fo n\in\NN\st \tr(T+2n)\in\cN'\}$, we have
a.s.
\bq
\liminf_{n\ri\iy}\frac{\ln(I_{T+2n})}{n}&\geq& \gamma\eq
Lemma \ref{ab} permits extending these results to the statement of Proposition 
\ref{stop} (up to replacement of $R_2$ by $bR_2/a$.
\par\me
Now the proof of Theorem \ref{th1} is clear.
By iteration, introduce two sequences $(S_n)_{n\in\NN}$ and $(T_n)_{n\in\NN}$ of stopping times for the triangle Markov chain: 
 start with $S_0=0$ and for any $n\in\NN$, if $S_n$ has been defined, take
\bq
T_n&\df&\inf\{m> S_n\st \ln(I(\tr(m)))>R_2\}\\
S_{n+1}&\df& \inf\{m> T_n\st \ln(I(\tr(m)))<R_1\}\eq
Of course, if for some $n\in\NN$, $S_n=\iy$ then for any $m\geq n$, $S_m=T_m=\iy$.
Conversely, via Proposition~\ref{div}, we see that if $S_n<\iy$, then a.s., $T_n<+\iy$, so
the events $\{S_n<\iy\}$ and $\{T_n<\iy\}$ are the same, up to a negligible set.
For $n\in\NN$, let us define the event
\bq
E_n&\df& \{S_n<\iy\hbox{ and } S_{n+1}=\iy\}\\
&=&\{T_n<\iy\hbox{ and }\fo m\in\NN,\, \tr(T_n+m)\in\cN'\}\eq
Up to conditioning on $\{S_n<\iy\}$, Lemma \ref{stop} shows that
\bq
\PP[S_{n+1}=\iy\vert S_n<\iy]\ =\ 
\PP[E_n\vert S_n<\iy]\ \geq\ 1/2\eq
thus it follows easily that $\PP[\cup_{n\in\NN}E_n]=1$.
Lemma \ref{stop} also shows that on all the $E_n$, the sequence $(I_m^{-1})_{m\in\NN}$
converges exponentially fast to zero with rate at least $\gamma$. Now the bound (\ref{yI})
implies  the validity of Theorem \ref{th1} with $\chi=\gamma/2$. 
\begin{rem}\label{In}
Let $\gamma_2\df F(1/2)=\min_{x\in[0,1/2]}F(x)$. A closer look at the proof of Proposition \ref{neigh}
shows that for any $\gamma<\gamma_2$, we can find a neighborhood $\cN$
of the set of flat triangles such that the lower bound of Proposition \ref{neigh} is satisfied.
By the above arguments, it follows that Theorem \ref{th1} also holds
with $\chi=\gamma_2/2$, so we win a factor $1/2$.\\
But one can go further. For $N\in\NN\setminus\{0,1\}$ and $x\in[0,1/2]$, consider 
\bq
F_N(x)&\df& \frac1{6^N}\sum_{(i_1,...,i_N)\in\lin 1,6\rin^N}
\ln\lt(G(i_N,z_{i_{N-1}}\circ \cdots\circ z_{i_1}(x))\cdots G(i_2,z_{i_1}(x))G(i_1,x)\rt)\\
&=&\EE_x\lt[\sum_{n\in\lin 0,N-1\rin}\ln(G(I_{n+1},Z_n))\rt]\eq
where $(\iota_n)_{n\in\NN^*}$ is a sequence of independent random variables uniformly distributed on $\lin 1,6\rin$
and $(Z_n)_{n\in\NN}$ is the Markov chain starting from $x$ ($Z_0\df x$) constructed from $(\iota_n)_{n\in\NN^*}$ through the relations
\bqn{ZI}
\fo n\in\NN,\qquad Z_{n+1}&\df& z_{\iota_{n+1}}(Z_n)\eqn
Then define
\bq
\gamma_N&\df&\min_{x\in[0,1/2]}F_N(x)\eq
An easy extension of the previous proof shows that
Theorem \ref{th1} holds
with $\chi=\gamma_N/N$
and consequently with $\chi=\lim_{N\ri\iy}\gamma_N/N$.
The quantity $\gamma_N/N$  converges due to the weak convergence of
the Markov chain $(Z_n)_{n\in\NN}$, uniformly in its
initial distribution, as we will show  in the next section.
Indeed, if $\mu$ is the attracting invariant probability associated to $(Z_n)_{n\in\NN}$,
we will see that for any $x\in[0,1/2]$,
\bq
\lim_{n\ri\iy}\EE_x\lt[\ln(G(\iota_{n+1},Z_n))\rt]&=&L\eq
with \bqn{L}
L&\df&\frac16\sum_{i\in\lin 1,6\rin}\int \ln(G(i,x))\, \mu(dx)\eqn
It follows from Cesaro's lemma that
\bq
\lim_{n\ri\iy}\frac{F_N(x)}N&=&\lim_{n\ri\iy}\frac1N\sum_{n\in\lin 0,N-1\rin}
\EE_x\lt[\ln(G(\iota_{n+1},Z_n))\rt]\\
&=&L
\eq
Since this convergence holds uniformly in $x\in[0,1/2]$,
we get that 
Theorem \ref{th1} is satisfied
with $\chi=L$.
In section 7, we will numerically evaluate that $L\approx 0.07$.
\end{rem}

\section{Ergodicity of the limit flat Markov chain}

This section studies  the limit flat Markov chain $Z\df(Z_n)_{n\in\NN}$. First we will see that it admits a unique invariant probability $\mu$ and that it converges exponentially fast to $\mu$ in the Wasserstein distance. Next we will show that $\mu$ is continuous and that its support is the whole state space $[0,1/2]$.\par\me
We begin by describing the kernel of $Z$ given in (\ref{Mz}) in the language of
iterated random functions.
\begin{lem}\label{homz}
With the notation of the previous sections, we have for all $x\in[0,1/2]$,
\bq
\begin{array}{rclcrcl}
 z_1(x)&=&\frac{3x}{2+2x}&\qquad&
 z_2(x)&=&\frac{3x}{2-x}\un_{x<2/7}+\frac{2-4x}{2-x}\un_{x\geq2/7}\\
 \noalign{\vskip-3mm}\\
 z_3(x)&=&\frac{1+x}{3-3x}\un_{x<1/5}+\frac{2-4x}{3-3x}\un_{x\geq1/5}&\qquad&
 z_4(x)&=&\frac{1+x}{4-2x}\\
  \noalign{\vskip-3mm}\\
 z_5(x)&=&\frac{1-2x}{4-2x}&\qquad&
 z_6(x)&=&\frac{1-2x}{3}
 \end{array}
 \eq
\end{lem}
\proof
These are immediate computations, based on the fact that for any flat triangle, the abscissa of the characteristic point is the ratio of the shortest edge by the longest edge.
For instance the lengths of the edges of the triangle $\tr_2$ are
${L_1}/2$, ${l_1}/3$ and ${2l_3}/{3}$ with 
$L_1 = x$, 
$l_1= 1-\frac{x}2$
and
$l_3= \frac12-x$, which leads to the above expression for $z_2(x)$.
\wwtbp
To see that the Markov kernel $M$ of $Z$ is ergodic, in the sense that it admits an invariant  and attracting  probability, we apply a result due to Barnsley and Elton \cite{MR932532}:
let $S$ be a compact segment of $\RR$ (more generally it can be a complete, separable metric space) on which we are given $n$ Lipschitz
functions $f_i\st S\ri S$, for $i\in\lin 1,n\rin$.
Let $p=(p_i)_{i\in\lin 1,n\rin}$ be a probability on $\lin 1,n\rin$
and consider the Markov kernel $N$ from $S$ to $S$ given by
\bqn{rfMc}
\fo x\in S,\qquad N(x,\cdot)&\df& \sum_{i\in\lin 1, n\rin}p_i\delta_{f_i(x)}\eqn
Then under the assumption that there exists a constant $r<0$ such that 
\bqn{BE2}
\fo x\not=y\in S,\qquad\sum_{i\in\lin 1, n\rin}p_i\ln\lt(\frac{\vert f_i(y)-f_i(x)\vert}
{\lve y-x\rve}\rt)&\leq& r\eqn
the kernel $N$ is ergodic: it admits a unique invariant and attracting probability $\mu$,  satisfying $\mu N=\mu$ and for any probability $\nu$ on $S$, $\lim_{n\ri\iy}\nu N^n=\mu$ (in the weak topology).
Furthermore, Barnsley and Elton \cite{MR932532} show that
there exists $q\in(0,1]$ and $\rho\in(0,1)$ such that
\bqn{BE}
\fo x,y\in S,\qquad \sum_{i\in\lin 1,n\rin}p_i\lve f_i(y)-f_i(x)\rve^q&\leq&\rho\lve y-x\rve^q\eqn
Let us rewrite this bound in a more probabilistic way.
Let $(\iota_n)_{n\in\NN}$ be a sequence of independent random variables taking values in $\lin 1,n\rin$
with distribution $(p_i)_{i\in\lin 1,n\rin}$.
For any $x\in S$, we denote by $U^x\df(U^x_n)_{n\in\NN}$ the stochastic chain constructed as follow: $U^x_0=x$
and for any $n\in\NN$, $U^x_{n+1}=f_{\iota_{n+1}}(U^x_n)$. This is a Markov chain with transition kernel $N$.
This construction enables us to couple together all the Markov chains $U^x$, for $x\in S$.
Then the above bound can be written
\bq
\fo x,y\in S,\qquad \EE[\lve U_1^y-U_1^x\rve^q]&\leq& \rho\lve y-x\rve^q\eq
and admits an immediate extension:
\bq
\fo n \in\NN,\,\fo x,y\in S,\qquad \EE[\lve U_n^y-U_n^x\rve^q]&\leq& \rho^n\lve y-x\rve^q\eq
This leads us to consider the Wasserstein distance $D$ between probability measures
on $S$: if $\nu_1$ and $\nu_2$ are two such measures,
\bq
D(\nu_1,\nu_2)&\df& \sup_{f\in\cL(1)}\lve \nu_1[f]-\nu_2[f]\rve\eq
where $\cL(1)$ is the set of Lipschitz functions on $S$ whose Lipschitz constant is less (or equal) than 1.
We may now show that
\begin{lem}\label{Wass}
Under the above assumption (\ref{BE2}), we have for any  $n\in\NN$ and any $x\in S$,
\bq
D(N^n(x,\cdot),\mu)&\leq& \diam(S)\rho^n\eq
where $\rho$ and $q$ are as in (\ref{BE}) and $\diam(S)$ is the diameter of $S$.
It follows that $U^x$ satisfies the law of large numbers:
for any continuous function $f$ on $S$,
we have a.s.,
\bqn{ff}
\lim_{N\ri\iy}\frac1{N+1}\sum_{n\in\lin 0,N\rin}f(U^x_n)&=&\mu[f]\eqn
\end{lem}
\proof
Let $f\in\cL(1)$. We compute that
\bq
\lve N^n(x,f)-\mu[f]\rve&=&\lve \int \mu(dy)(N^n(x,f)-N^n(y,f))\rve\\
&\leq&\sup_{x,y\in S}\lve N^n(x,f)-N^n(y,f)\rve\\
&=&\sup_{x,y\in S}\lve \EE[f(U^x_n)-f(U^y_n)]\rve\\
&\leq&\sup_{x,y\in S}\EE[\lve U^x_n-U^y_n\rve]\\
&=&\diam(S)\sup_{x,y\in S}\EE\lt[\frac{\lve U^x_n-U^y_n\rve}{\diam(S)}\rt]\\
&\leq& \diam(S)\sup_{x,y\in S}\EE\lt[\lve \frac{U^x_n-U^y_n}{\diam(S)}\rve^q\rt]\\
&\leq& \diam(S)^{1-q}\sup_{x,y\in S}\rho^n\lve y-x\rve^q\\
&\leq& \diam(S)\rho^n\eq
The announced bound follows by taking the supremum over all functions $f\in\cL(1)$.
The law of large numbers is deduced from a traditional martingale argument based on the existence
of a bounded solution to the Poisson equation. More precisely,
for $f\in\cL(1)$, we can define
\bq
\fo x\in S,\qquad \varphi(x)&\df& \sum_{n\in\NN}\EE[f(U^x_n)-\mu[f]]\eq
since the r.h.s.\ converges exponentially fast and uniformly with respect to $x\in S$.
Furthermore we easily see that $\varphi$ is a Lipschitz function and that
it is a solution to the Poisson equation
\bq
\lt\{\begin{array}{rcl}
\fo x\in S,\qquad \varphi(x)-N(x,\varphi)&=&f(x)-\mu[f]\\
\mu[\varphi]&=&0
\end{array}\rt.\eq
This enables us to write for any $n\in\NN$,
\bq
f(X_0)+f(X_1)+\cdots+f(X_n)&=&(n+1)\mu[f]+\varphi(X_0)-\varphi(X_{n+1})+\cM_{n+1}\eq
where $(\cM_n)_{n\in\NN}$ is a martingale whose increments are bounded. 
The law of large numbers for functions $f$ belonging to $\cL(1)$ then follows from the well-known fact that 
$\cM_n/n$ converges a.s.\ to zero. It is also true for all Lipschitz functions $f$.
Next,  given a continuous function $f$ on $S$ and $m \in\NN^*$, by usual approximations, it is possible to find
a Lipschitz function $\wi f_m$  on $S$ such that $\lVe f-\wi f_m\rVe_{S,\iy}\leq1/m$,
where $\lVe \cdot\rVe_{S,\iy}$ is the uniform norm on $S$.
It follows that on a mesurable set $\Omega_m$ of probability 1,
 \bq
\limsup_{N\ri\iy}\frac1{N+1}\sum_{n\in\lin 0,N\rin}f(U^x_n)&\leq&\mu[f]+2/m\\
\liminf_{N\ri\iy}\frac1{N+1}\sum_{n\in\lin 0,N\rin}f(U^x_n)&\geq&\mu[f]-2/m\eq
Thus on the set $\cap_{m\in \NN^*}\Omega_m$ of full probability, (\ref{ff}) is true.
\wwtbp
Let us discuss condition (\ref{BE2}).
Note that since the functions $f_i$, for $i\in\lin1,n\rin$, are Lipschitz, they are
 absolutely continuous. Let us write
$f'_i$ for their respective weak derivatives. 
By letting $y$ and $x$ become close in  criterion (\ref{BE2}),
we get that almost everywhere in $x\in S$, 
\bqn{cn}\sum_{i\in\lin 1, n\rin}p_i\ln(\vert f'_i(x)\vert)&\leq& r\eqn
Condition (\ref{cn}) is not sufficient to insure that the kernel $N$ is ergodic.
Consider the following example with $S=[0,1]$, 
$n=2$ and the functions $f_1$ and $f_2$ defined by
\bq
\fo x\in[0,1],\qquad f_i(x)&\df& \lt\{\begin{array}{ll}\min(2x,1)&\hbox{, if $i=1$}\\
\max(0,-1+2x)\hbox{, if $i=2$}\end{array}\rt.\eq
In this case (\ref{cn}) is even satisfied with $r=-\iy$ and the set of invariant probability measures is $\{a\delta_0+(1-a)\delta_1\st a\in[0,1]\}$,
so none of them can be attractive (but the law of a corresponding Markov chain converges exponentially fast
to one of the invariant probability measures).\par\sm
Nevertheless,
under some circumstances, the necessary condition (\ref{cn}) is also sufficient.
This is  the case if for all the functions $\lve f_i'\rve$, with $i\in\lin1,n\rin$,  there exist $a_i,b_i\in\RR$ such that almost everywhere (a.e.) in $x\in S$, 
\bqn{cs}
\lve f_i'(x)\rve&=&(a_ix+b_i)^{-2}\eqn
(in particular $-b_i/a_i$ cannot belong to $S$ otherwise $f_i'$ would not be integrable over this interval).
Indeed in this situation we can write that for any $x<y\in S$,
\bq
\frac{\vert f_i(y)-f_i(x)\vert}
{\lve y-x\rve}&=&
\frac{1}{\lve y-x\rve}\lve \int_x^y f'_i(z)\,dz\rve\\
&\leq& \frac{1}{\lve y-x\rve} \int_x^y \lve f'_i(z)\rve\,dz\\
&=&\frac{1}{\lve y-x\rve} \int_x^y \frac1{(a_iz+b)^2}\,dz\\
&=&\frac{1}{a_i\lve y-x\rve}\lt(\frac1{a_ix+b_i}-\frac1{a_iy+b_i}\rt)\\
&=&\frac{1}{a_i\lve y-x\rve}\frac{a_i(y-x)}{(a_iy+b_i)(a_ix+b_i)}\\
&=&\frac1{\lve a_iy+b_i\rve\lve a_ix+b_i\rve}\\
&=&\sqrt{\lve f_i'(y)\rve \lve f_i'(x)\rve}
\eq
where the last equality has to be understood a.e.
It follows that, at least for a.e.\  $x,y\in S$,
\bqn{lnln}
\ln\lt( \frac{\vert f_i(y)-f_i(x)\vert}
{\lve y-x\rve}\rt)&\leq& \frac{\ln(\lve f_i'(y)\rve)+\ln(\lve f_i'(x)\rve)}{2}\eqn
and consequently
\bq
\sum_{i\in\lin 1,n\rin}p_i\ln\lt( \frac{\vert f_i(y)-f_i(x)\vert}
{\lve y-x\rve}\rt)&\leq& \frac12\sum_{i\in\lin 1,n\rin}p_i\ln(\lve f_i'(y)\rve)+\frac12\sum_{i\in\lin 1,n\rin}p_i\ln(\lve f_i'(x)\rve)\eq
A formula which enables passing from (\ref{cn}) to (\ref{BE2}). It  only has to be checked for a.e.\ $x,y\in S$.
\par\me
It is time now to come back to the flat triangle Markov chain. Consider the setting where $N=M$, i.e.\ $S=[0,1/2]$,
$n=6$, $f_i=z_i$ for $i\in\lin 1, 6\rin$ and $p$ the uniform distribution on $\lin 1, 6\rin$.
Now condition (\ref{cs}) is satisfied.
Since $z_2'(2/7-)=-z_2'(2/7+)$ and $z_3'(1/5-)=-z_3'(1/5+)$, we see that
$\vert z'_i\vert(x)$ can be defined everywhere (a convention that we will adopt from now on).
Indeed, we compute that for any $x\in[0,1/2]$,
\renewcommand{\arraystretch}{2.2}
\bq\begin{array}{rclcrcl}
 \lve z'_1\rve(x)&=&\di\frac{3}{2(1+x)^2},&\qquad&
 \lve z'_2\rve(x)&=&\di\frac{6}{(2-x)^2}\cr
 \lve z'_3\rve(x)&=&\di\frac{2}{3(1-x)^2},&\qquad&
 \lve z'_4\rve(x)&=&\di\frac{3}{2(2-x)^2}\cr
 \lve z'_5\rve(x)&=&\di\frac{3}{2(2-x)^2},&\qquad&
 \lve z'_6\rve(x)&=&\di\frac{2}{3}\end{array}\eq
 \renewcommand{\arraystretch}{1}
Unfortunately (\ref{cn}) is not true and surprisingly it is a computation we have already encountered:
comparing with (\ref{G}), we see that 
 \bq
 \fo i\in\lin 1,6\rin,\,\fo x\in[0,1/2],\qquad
 \lve z'_i\rve(x)&=&\frac1{G^2(i,x)}\eq
thus by the observation before Proposition \ref{neigh}, we know
that $\sum_{i\in\lin 1, 6\rin}\ln(\vert z'_i(x)\vert)$ is positive for $x$ near $1/2$.
As in section 3, we get around this difficulty by iterating the kernel $M$ one more time (this trick was also used by Barnsley and Elton in one example of their paper \cite{MR932532}). So we consider
$N=M^2$,
 namely $S=[0,1/2]$,
$n=36$, $f_{i,j}=z_i\circ z_j$ for $(i,j)\in\lin 1, 6\rin^2$ and $p$ the uniform distribution on $\lin 1, 6\rin^2$.
The advantage is that we have 
for any $i,j\in\lin 1,6\rin$ and any $x\in [0,1/2]$,
 \bq 
 \lve f'_{i,j}\rve(x)&=&\lve z'_i\rve (z_j(x))\lve z'_i\rve(x)\\
 &=&\lt(G(i,z_j(x))G(j,x)\rt)^{-2}\eq
Thus 
\bq
\fo x\in[0,1/2],\qquad \sum_{i,j\in\lin 1,6\rin}\ln(\vert f'_{i,j}\vert(x))&=&-2F(x)\eq
and in particular the left hand side is negative due to Proposition \ref{neigh} (it is even increasing as a function of $x\in[0,1/2]$ according to the observation made at the end of section 3).
But (\ref{cs}) is no longer satisfied by the functions $f_{i,j}$.
To avoid this problem, we come back directly to the bound (\ref{lnln}): for $i,j\in\lin 1,6\rin$
and $y>x\in[0,1/2]$, we write
\bq
\ln\lt( \frac{\vert f_{i,j}(y)-f_{i,j}(x)\vert}
{\lve y-x\rve}\rt)&=& 
\ln\lt( \frac{\vert z_{i}(z_j(y))-z_{i}(z_j(x))\vert}
{\lve z_j(y)-z_j(x)\rve}\rt)+\ln\lt( \frac{\vert z_{j}(y)-z_{j}(x)\vert}
{\lve y-x\rve}\rt)\\
&\leq& \frac{\ln(\lve z_i'(z_j(y))\rve)+\ln(\lve z_i'(z_j(x))\rve)}{2}
+\frac{\ln(\vert z_j'(y)\vert)+\ln(\vert z_j'(x)\vert)}{2}\\
&=&\frac{\ln(\vert f_{i,j}'(y)\vert)+\ln(\vert f_{i,j}'(x)\vert)}{2}
\eq
In this situation we can also come back from (\ref{cs}) to (\ref{cn})
and the results of Barnsley and Elton \cite{MR932532} insure that the iterated Markov kernel $M^2$
 is ergodic. 
To come back from $M^2$ to $M$ is not difficult:
\begin{pro}\label{slln}
The kernel $M$ is ergodic and the Markov chain $Z$ satisfies the strong law of large numbers.
\end{pro}
\proof
Let $\mu$ be the attracting and invariant probability for $M^2$.
Then we have $(\mu M)M^2=(\mu M^2)M=\mu M$, so $\mu M$ is invariant for $M^2$
and by uniqueness it follows that $\mu M=\mu$. Next for any probability measure $\nu$ on $[0,1/2]$,
the (weak) limit set of $(\nu M^n)_{n\in\NN}$ is included in $\{\mu, \mu M\}=\{\mu\}$,
so $\mu$ is also attracting for $M$ and the uniqueness of $\mu$ as the invariant probability of $M$ follows.
Finally the strong law of large numbers for $Z$ can be deduced from that
of the two Markov chains $(Z_{2n})_{n\in\NN}$ and $(Z_{1+2n})_{n\in\NN}$.
\wwtbp
There was a cruder way to deduce the ergodicity of $M$:
 \begin{rem}\label{LipDF}
Diaconis and Freedman  \cite{MR1669737} consider a simpler criterion 
for ergodicity of a random function Markov kernel (\ref{rfMc}): for $i\in\lin 1, n\rin$, let
$K_i\df\sup_{x\not=y}{\lve f_i(y)-f_i(x)\rve}/{\lve y-x\rve}$ be
 the Lipschitz constant of $f_i$ and
assume that there exists a constant $r<0$ such that 
\bqn{DF}
\fo x,y\in S,\qquad\sum_{i\in\lin 1, n\rin}p_i\ln\lt(K_i\rt)&\leq& r\eqn
then the kernel $N$ is ergodic and Diaconis and Freedman  \cite{MR1669737} show that the convergence
is exponentially fast in the Prokhorov distance (but for us the Wasserstein distance is more convenient because in the end we would like to couple the two Markov chains $(X_n,Y_n)_{n\in\NN}$
 and $(Z_n)_{n\in\NN}$).
Of course condition (\ref{DF}) implies (\ref{BE2}). Since for $i\in\lin 1,n\rin$, $K_i$ is the essential supremum
of $\lve f'_i(x)\rve$, (\ref{DF}) corresponds to the exchange of essential supremum and sum
in (\ref{BE2}). Let us now come back to our flat triangle Markov chain.
From the previous considerations, (\ref{DF}) cannot be satisfied with $N=M$.
It does not work either with $N=M^2$, so this is an example were the criterion (\ref{BE2}) is fulfilled
while (\ref{DF}) is not.
But condition (\ref{DF}) is satisfied with $N=M^3$, namely 
$S=[0,1/2]$,
$n=216$, $f_{i,j,k}=z_i\circ z_j\circ z_k$ for $(i,j,k)\in\lin 1, 6\rin^3$ and $p$ the uniform distribution on $\lin 1, 6\rin^3$. For the details of the underlying numerical computations, we refer to section 7.
 \end{rem}
To finish  we prove two properties of $\mu$ which will be needed in next section. The first one is
\begin{lem}\label{zero}
The probability $\mu$ contains no atom, in particular
$\mu(\{0\})=0$.
\end{lem}
\proof
The proof needs a few steps and notation. Let us define 
\bq\mu^*&\df&\sup\{\mu(\{x\})\st x\in[0,1/2]\}\\
S^*&\df& \{x\in[0,1/2]\st \mu(\{x\})=\mu^*\}\eq
 and for any $x\in[0,1/2]$,
\bq
\bar{S}(x)&\df&\{(i,y)\in\lin 1,6\rin\times [0,1/2]\st z_i(y)=x\}\\
S(x)&=&\{y\in[0,1/2]\st \ex i\in\lin 1,6\rin \hbox{ with } z_i(y)=x\}\eq
The first step is:
\par\sm
$\bullet$ We have $\mu(\{0\}=\mu(\{1/2\})\leq \mu^*/2$.\par\sm
By invariance of $\mu$ we can write that
\bq
\mu(\{0\})&=&\mu(M[\un_{\{0\}}])\\
&=&\frac16\sum_{(i,y)\in\bar{S}(0)}\mu(\{y\})\\
&=&\frac26\mu(\{0\})+\frac46\mu(\{1/2\})\eq
and this relation implies that $\mu(\{0\})=\mu(\{1/2\})$.
Next consider the point $1/2$. We get
\bq
\mu(\{1/2\})&=&\frac16\sum_{(i,y)\in\bar{S}(1/2)}\mu(\{y\})\\
&=&\frac26\mu(\{1/2\})+\frac16\mu(\{1/5\})+\frac16\mu(\{2/7\})\\
&\leq& \frac26\mu(\{1/2\})+\frac13\mu^*\eq
so it follows that $\mu(\{1/2\})\leq \mu^*/2$.\par\sm
The next step is:\par\sm
$\bullet$ For any $x\in S^*$, we have $S(x)\subset S^*\cup\{0\}$.
\par\sm
Looking at the graphs of the functions $z_i$, for $i\in\lin 1,6\rin$ (see Figure 1 in section 7), we get that
\bq
\fo x\in[0,1/2],\qquad \card(\bar{S}(x))&=&\lt\{\begin{array}{ll}
4&\hbox{, if $x=1/2$}\\
7&\hbox{, if $x=1/4$ or $x=1/3$}\\
6&\hbox{, otherwise}\end{array}\rt.\eq
So consider $x\in S^*\setminus\{1/4,1/3,1/2\}$, writing 
\bq
\mu(\{x\})&=&\frac16\sum_{(i,y)\in \bar{S}(x)}\mu(\{y\})\\
&\leq& \mu^*\eq
it appears that equality is possible only if $\mu(\{y\})=\mu^*$
for all $y\in S(x)$, namely $S(x)\subset S^*$.\\
We study now the three particular cases of $1/4$, $1/3$ and $1/2$.\\
- For $1/2$: as seen in the first step, $1/2\in S^*$ implies that $\mu^*=0$,
so $S^*=[0,1/2]$ and the inclusion $S(1/2)\subset S^*$ is trivial.\\
- For $1/4$: there exist five distinct points $y'_1,y'_2,y'_3,y'_4,y'_5\in[0,1/2]$
such that we have 
\bq
S(1/4)&=&\{(4,0),(5,0),(1,y'_1),(2,y'_2),(2,y'_3),(3,y'_4),(6,y'_5)\}\eq
so by invariance of $\mu$ we get
\bq
\mu(\{1/4\})&=&\frac16\lt(2\mu(\{0\})+\sum_{i\in\lin 1,5\rin}\mu(\{y'_i\})\rt)
\\
&\leq& \frac16(2\mu(\{0\})+5\mu^*)\eq
and since we know that $\mu(\{0\})\leq \mu^*/2$, the equality $\mu(\{1/4\})=\mu^*$
is possible only if $\mu(\{y'_i\})=\mu^*$ for $i\in\lin 1,5\rin$, so we can conclude that $S(1/4)\subset S^*\cup\{0\}$.\\
- The same argument holds for $1/3$ (even if $1/3\in S(1/3)$).\par
For the last step, let us denote by $\wi z_2$ the restriction
of $z_2$ to $[0,2/7]$. This mapping is one-to-one from $[0,2/7]$ to $[0,1/2]$
and we denote  its inverse by $\wi z_2^{-1}$.\par\sm
$\bullet$ For $x\in(0,1/2]$, the set $\{\wi z_2^{-n}(x)\st n\in\NN\}$ is infinite, so $S^*$ is infinite.\par\sm
The first assertion comes from the fact that for any $x\in (0,1/2]$,
$0<\wi z_2^{-1}(x)<x$, so $(\wi z_2^{-n}(x))_{n\in\NN}$ is indeed a decreasing sequence (converging to 0).
By the first step, $S^*$ cannot be reduced to $\{0\}$, so there exists $x\in S^*\setminus\{0\}$.
By the previous step, the sequence $(\wi z_2^{-n}(x))_{n\in\NN}$ is included in $S^*$, since
none of its elements can be equal to $0$.
It follows that $S^*$ is infinite.\par\sm
Of course the last statement implies that $\mu^*=0$, because $\mu$ is a probability measure.
\wwtbp
If all of the functions $z_i$, for $i\in\lin 1,6\rin$, were strictly monotone, the fact that $\mu(\{0\})=0$ could have been
deduced more directly from the uniqueness of $\mu$ and Theorem 2.10 from Dubins and Freedman \cite{MR0193668}.
The second piece of information we  need about $\mu$ is a direct consequence of a result of the latter paper.
\begin{lem}\label{support}
The support of $\mu$ is the whole segment $[0,1/2]$.
\end{lem}
\proof
By Theorem 4.9 of Dubins and Freedman \cite{MR0193668}, the support of $\mu$ is the whole segment $[0,1/2]$
if we can cover it by the images of the functions $z_i$ which are strict contractions.
But this is the case here, since $z_4$ and $z_5$ are strict contractions and $z_4([0,1/2])=[1/4,1/2]$
and $z_5([0,1/2])=[0,1/4]$.\wwtbp

\section{More on the asymptotic behaviour}

Our main goal here is to prove Theorems \ref{th2} and \ref{th3}. 
The underlying tool is to couple the Markov chains $(X_n,Y_n)_{n\in\NN}$ and $(Z_n)_{n\in\NN}$
to take advantage of the information we have on the chain $(Z_n)_{n\in\NN}$.\par\me
A natural coupling between the above chains is based on the construction alluded to
in Remark \ref{In}. 
Assume that $(X_0,Y_0)$ and $Z_0$ are given and let $(\iota_n)_{n\in\NN^*}$ be a sequence of independent random variables uniformly distributed on $\lin 1,6\rin$ and independent from the previous initial conditions.
We consider $(Z_n)_{n\in\NN}$ constructed as in (\ref{ZI}) and similarly we iteratively define 
$(X_n,Y_n)_{n\in\NN}$ via
\bq
\fo n\in\NN,\qquad (X_{n+1},Y_{n+1})&\df& (x_{\iota_{n+1}}(X_n),y_{\iota_{n+1}}(Y_n))\eq
In these relations, the indices refer to the conventions made in (\ref{triangles}) and (\ref{zx}).
A first simple property of this coupling is
\begin{lem}\label{mauvais}
The random variables $\lve X_n-Z_n\rve$ converge in probability to zero as $n$ goes to infinity:
\bq
\fo \epsilon >0,\qquad\lim_{n\ri\iy}\PP[\lve X_n-Z_n\rve>\epsilon]&=&0\eq
\end{lem}
\proof
First we iterate Lemma \ref{QM} to show that with $K'\df K^2+8K/3>0$, we have
\bq
\fo i,j\in\lin 1,6\rin,\,\fo (x,y)\in\cD,\qquad 
\lve x_i(x_j(x,y),y_j(x,y))-z_i(z_j(x))\rve&\leq& K'{y}
\eq
Indeed, taking into account that all the functions $z_i$, $i\in\lin 1,6\rin$ have a Lipschitz constant
less than (or equal to) 8/3, we deduce that for any $i,j\in\lin 1,6\rin$ and any $(x,y)\in\cD$,
\bq
\lefteqn{\lve x_i(x_j(x,y),y_j(x,y))-z_i(z_j(x))\rve}\\
&\leq&\lve x_i(x_j(x,y),y_j(x,y))-z_i(x_j(x,y))\rve
+\lve z_i(x_j(x,y))-z_i(z_j(x))\rve\\
&\leq& K\lve y_j(x,y)\rve+\frac83\lve x_j(x,y)-z_j(x)\rve\\
&\leq&K^2{y}+\frac{8K}3 y\eq
Let $q\in(0,1]$ and $\rho\in(0,1)$ as in (\ref{BE}) but relative to the kernel $N=M^2$.
Then for any $n\in\NN$, we can write
\bq
\lefteqn{\EE[\vert X_{n+2}-Z_{n+2}\vert^q\vert X_n,Y_n,Z_n]}\\
&=&\EE[
\vert x_{\iota_{n+2}}(x_{\iota_{n+1}}(X_n,Y_n),y_{\iota_{n+1}}(X_n,Y_n))-z_{\iota_{n+2}}(z_{\iota_{n+1}}(Z_n))^q\vert \vert X_n,Y_n,Z_n]\\
&\leq& \EE[\vert
x_{\iota_{n+2}}(x_{\iota_{n+1}}(X_n,Y_n),y_{\iota_{n+1}}(X_n,Y_n))-z_{\iota_{n+2}}(z_{\iota_{n+1}}(X_n))\vert^q \vert X_n,Y_n,Z_n]
\\&&+\EE[\vert z_{\iota_{n+2}}(z_{\iota_{n+1}}(X_n))-z_{\iota_{n+2}}(z_{\iota_{n+1}}(Z_n))\vert^q \vert X_n,Y_n,Z_n]\\
&\leq&(K')^qY_n^{{q}}+\rho \vert X_n-Z_n\vert^q\eq
For $n\in\NN$, denote 
\bq
a_n&\df& \EE[\vert X_n-Z_n\vert^q]\\
b_n&\df&(K')^q\EE[Y_n^{{q}}]\eq
After integration, the above bound leads to
\bq
\fo n\in\NN,\qquad a_{n+2}&\leq& \rho a_n+b_n\eq
We deduce that
\bqn{ab2}
\fo n\in\NN,\qquad a_{2n}&\leq &a_0\rho^n+\sum_{m\in\lin 0, n-1\rin}b_{2m}\rho^{n-1-m}\eqn
where  $\lim_{n\ri\iy} a_{2n}=0$ is a consequence of
$\lim_{n\ri\iy} b_{n}=0$. A similar computation shows that this latter condition
also implies that $\lim_{n\ri\iy} a_{2n+1}=0$, i.e.\ at the end we will be insured of
\bq\lim_{n\ri\iy} \EE[\vert X_n-Z_n\vert^q]&=&0\eq
and thus of the announced convergence in probability.\\
But we already know that $(Y_n)_{n\in\NN}$ converges a.s.\ to zero
and since this sequence is uniformly bounded, we see by the dominated convergence theorem
that $\lim_{n\ri\iy} b_{n}=0$. \wwtbp
Now Theorem \ref{th2} follows quite easily:
\par\me\noindent\textbf{Proof of Theorem \ref{th2}}\par\sm\noindent
For $n\in\NN$,  denote by $A_n'$ (respectively $A_n''$) the angle between $[(0,0),(X_n,Y_n)]$ and $[(X_n,Y_n),(X_n,0)]$
(resp.\ $[(X_n,0),(X_n,Y_n)]$ and $[(X_n,Y_n),(1,0)]$), so that $A_n=A_n'+A_n''$.
Since the length of $[(X_n,0),(1,0)]$ is larger than 1/2 and $Y_n$ converges a.s.\ to zero for
large $n\in\NN$, it is clear that $A_n''$ converges a.s.\ to $\pi/2$.
Furthermore, we have $\tan(A_n')=X_n/Y_n$, so to see that 
$A_n$ converges in probability toward $\pi$, we must
see that $Y_n/X_n$ converges in probability toward 0.
Let $\epsilon,\eta >0$ be given, we have for any $n\in\NN$,
\bqn{PPPP}
\nonumber\PP[Y_n/X_n\geq\epsilon]&=&\PP[Y_n/X_n\geq\epsilon,X_n> 2\eta]+\PP[Y_n/X_n\geq\epsilon,X_n\leq 2\eta]\\
\nonumber&\leq&\PP[Y_n\geq2\epsilon\eta]+\PP[X_n\leq 2\eta]\\
&\leq& \PP[Y_n\geq2\epsilon\eta]+\PP[\vert X_n-Z_n\vert \geq\eta]+\PP[Z_n\leq\eta]\eqn
By letting $n$ going to infinity, taking into account that the stationary distribution $\mu$ of $Z$ is continuous,
we get 
\bq
 \limsup_{n\ri\iy}\PP[Y_n/X_n\geq\epsilon]&\leq& \lim_{n\ri\iy}\PP[Z_n\leq\eta]\\
&=& \mu([0,\eta])\eq
because as $n$ goes to infinity, the law of $Z_n$ converges weakly to $\mu$
and this probability gives weight 0 to the boundary $\{\eta\}$ of $(-\iy,\eta]$.
 Using again Lemma \ref{zero} and letting $\eta$ go to zero,
 we obtain that $\lim_{n\ri\iy}\PP[Y_n/X_n\geq\epsilon]=0$
 and consequently the announced convergence in probability.\wwtbp
 \begin{rem}\label{expoa}
 We don't know if $(A_n)_{n\in\NN}$ converges to $\pi$ a.s.
One way to deduce this result, via the Borel-Cantelli lemma, would be to show that for any given $\epsilon >0$,
\bqn{sumfi}
\sum_{n\in\NN}\PP[Y_n/X_n\geq\epsilon]&<&+\iy\eqn
In view of the above arguments, one of the main problems is that we have no bound
on the way $\mu([0,\eta])$ goes to zero as $\eta$ goes to zero. We would like to find $\alpha>0$
such that $\limsup_{\eta\ri 0_+}\mu([0,\eta])/\eta^\alpha$
\linebreak 
$<+\iy$, but we were not able to
prove such an estimate. 
If we knew that $\mu$ is absolutely continuous, Figure 7 in section 7 would suggest that this property
holds with $\alpha=1$ (and $\lim_{\eta\ri 0_+}\mu([0,\eta])/\eta\leq1$).
\end{rem}
In order to prove Theorem \ref{th3}, we  need two technical results.
In all that follows, let fix some $a\in[0,1/2]$ and $\epsilon>0$ and define
$\cO\df[a-\epsilon,a+\epsilon]\cap [0,1/2]$.
\begin{lem}\label{cO}
There exist $\eta >0$ and $N\in\NN^*$
such that 
\bq \inf_{z\in [0,1/2]}\PP_z[Z_N\in \cO]&\geq&\eta\eq
(the index $z$ means that $Z_0=z$).
\end{lem}
\proof
This is a consequence of Lemma \ref{Wass} applied to $M^2$:
there exists $\rho\in(0,1)$ such that
for any $z\in[0,1/2]$ and $n\in\NN$, $D(M^{n}(z,\cdot),\mu)\leq\rho^{\lfloor n/2\rfloor}/2$.
Let $\varphi$ be the function vanishing outside $(a-\epsilon,a+\epsilon)$, affine on $[a-\epsilon,a]$
and $[a,a+\epsilon]$ such that $\varphi(a)=\epsilon$.
By definition of $D$, we have
\bq
\fo z\in[0,1/2],\,\fo n\in\NN,\qquad \lve M^n(z,\varphi)-\mu[\varphi]\rve&\leq& \frac{ \rho^{\lfloor n/2\rfloor}}{2}\eq
 Since the support of $\mu$ is [0,1/2], we have
$\eta\df\mu[\varphi]>0$. So if we choose $N\in\NN$ such that
$\rho^{\lfloor N/2\rfloor}<\eta$, we get that for any $z\in[0,1/2]$, $\PP_z[X_{N}\in \cO]\geq
M^{N}(z,\varphi)>\eta/2$.\wwtbp
For the second technical result, we need further notation:
for $(x,y)\in\cD$ and $(i_1,i_2,\cdots, i_N)\in \lin 1,6\rin^N$, we denote $x_{i_1,i_2,\cdots, i_N}(x,y)$,
 $y_{i_1,i_2,\cdots, i_N}(x,y)$ and $z_{i_1,i_2,\cdots, i_N}(x)$ the values of $X_n$, $Y_n$ and $Z_n$ when $(X_0,Y_0,Z_0)=(x,y,z)$
and $(\iota_1,\iota_2,\cdots, \iota_N)=(i_1,i_2,\cdots, i_N)$. 
\begin{lem}\label{xN}
There exists a constant $K_N$ such that \bq
\fo (i_1,i_2,\cdots, i_N)\in \lin 1,6\rin^N,\,\fo (x,y)\in\cD,\qquad \lt\{
\begin{array}{rcl}
\lve x_{i_1,i_2,\cdots, i_N}(x,y)-z_{i_1,i_2,\cdots, i_N}(x)\rve&\leq& K_Ny\\
\lve y_{i_1,i_2,\cdots, i_N}(x,y)\rve&\leq& K_Ny\end{array}\rt.
\eq
\end{lem}
\proof
For $N=1$ this is just Lemma \ref{QM}. The general case is proven by an easy iteration, similar to the one already used in the proof of Lemma \ref{mauvais},
starting with
\bq
\lefteqn{\lve x_{i_1,i_2,\cdots, i_N}(x,y)-z_{i_1,i_2,\cdots, i_N}(x)\rve}\\
&\leq& \lve x_{i_N}(x_{i_1,i_2,\cdots, i_{N-1}}(x,y),y_{i_1,i_2,\cdots, i_{N-1}}(x,y)
-z_{i_N}(x_{i_1,i_2,\cdots, i_{N-1}}(x,y))\rve\\&&
+\lve z_{i_N}(x_{i_1,i_2,\cdots, i_{N-1}}(x,y))-z_{i_N}(z_{i_1,i_2,\cdots, i_{N-1}}(x))\rve\eq
\wwtbp
We  can now come to the
\par\me\noindent\textbf{Proof of Theorem \ref{th3}}\par\sm\noindent
Let $\cO'\df[a-2\epsilon,a+2\epsilon]\cap [0,1/2]$.
We want to show an analogous result to Lemma \ref{cO} but for the chain $(X_n,Y_n)_{n\in\NN}$,
namely to find
$\eta' >0$ and $N'\in\NN^*$
such that 
\bqn{cOO} \inf_{(x,y)\in\cD}\PP_{(x,y)}[X_{N'}\in \cO']&\geq&\eta'\eqn
(let us recall that under $\PP_{(x,y)}$, $(X_0,Y_0)=(x,y)$).
To do so, we first consider
$\eta$  and $N$ as in Lemma \ref{cO} and consider $\delta>0$ sufficiently small
such that $K_N\delta^{1/2^{N-1}}<\epsilon$. Then according to Lemmas \ref{cO} and \ref{xN},
we have 
 \bq \inf_{(x,y)\in\cD\st y<\delta}\PP_{(x,y)}[X_{N}\in \cO']&\geq&\eta\eq
 To extend this estimate to the whole domain $\cD$, we come back to (\ref{yI}) and (\ref{plonge})
 which enables us to find $N''\in\NN$ such that
 \bq
\eta''&\df&\inf_{(x,y)\in\cD}\PP_{(x,y)}[Y_{N''}<\delta]\ >\ 0\eq
Now the Markov property implies (\ref{cOO}) with $\eta=\eta'\eta''$ and $N'=N+N''$.\par
\sm
 Since this bound is uniform over $(x,y)\in\cD$, the sequence $(\un_{\cO'}(X_{nN'}))_{n\in\NN}$
 is stochastically bounded below by an independent family of Bernoulli variables of parameter
 $\eta'$ and we deduce that a.s.
 \bq
 \limsup_{n\ri\iy}\un_{\cO'}(X_n)&=&1\eq
 The announced result follows because $a\in[0,1/2]$ and $\epsilon >0$ are arbitrary.\wwtbp
 The details of the above proof are necessary because in general
 one cannot deduce from the convergence in probability of $\lve X_n-Z_n\rve$ to
 zero as $n$ goes to infinity the a.s.\ equality of the limit sets of $(X_n)_{n\in\NN}$
 and of $(Z_n)_{n\in\NN}$. This property rather requires the a.s.\ convergence of $\lve X_n-Z_n\rve$ to
 zero and this leads to the following observations:
 \begin{rem}
Coming back to  Remark \ref{expoa}, to prove (\ref{sumfi}) via (\ref{PPPP}),
 we are also missing an estimate of the kind
\bqn{B}
\ex K,p,\chi>0\st\fo n\in\NN,\qquad \EE[Y_n^p]&\leq& K\exp(-\chi n)\eqn
Blackwell \cite{Blackwell_barycentric} succeeded in obtaining such a bound (with $p=1/2$)
by exhibiting an appropriate supermartingale with the help of the computer,
see also the survey by Butler and Graham 
\cite{Butler_triangle}.
His result can be seen to imply Theorem \ref{th1}, with $\chi =0.04$. 
\\
Furthermore it allows for a more direct proof of Theorem \ref{th3}.
Indeed, 
if (\ref{B}) is satisfied for some $p>0$, then for any $q>0$,
\bq
\sum_{n\in\NN}\EE[Y_n^q]&<&+\iy\eq
(this is immediate for $q=p$ and use the Hölder inequality for $0<q< p$ and the elementary bound
$y^q\leq (\sqrt{3}/2)^{q-p}y^p$ for $y\in[0,\sqrt{3}/2]$ and $q>p$).
The arguments of the proof of Lemma \ref{mauvais} (especially (\ref{ab2}) and a similar relation for odd integers) then show
that
\bq
\sum_{n\in\NN}\EE[\vert X_n-Z_n\vert^q]&<+\iy\eq
and consequently that
$\vert X_n-Z_n\vert$ converges a.s.\ to zero. This is sufficient to deduce that 
almost surely the limit set of $(X_n)_{n\in\NN}$ coincides with that of $(Z_n)_{n\in\NN}$,
thus
the law of large numbers for 
$Z$ and Lemma \ref{support} imply Theorem \ref{th3}.\par
In the same spirit, one can  go further toward justifying the assertion made in the introduction that
asymptotically $(X_n)_{n\in\NN}$ is almost Markovian.
Let us introduce the supremum distance $S$ on $[0,1/2]^\NN$, seen as the set of trajectories from $\NN$ to $[0,1/2]$:
\bq
\fo x\df(x_n)_{n\in\NN},\,z\df(z_n)_{n\in\NN}\in[0,1/2]^\NN,\qquad S(x,z)&\df& \sup_{n\in\NN}\lve x_n-z_n\rve\eq
For $m\in\NN$, let $X_{\lin m,\iy\lin}=(X_{m+n})_{n\in\NN}\in[0,1/2]^\NN$  and consider
\bq
s_m&\df& \inf\EE[S(X_{\lin m,\iy\lin},Z)]\eq
where the infimum is taken over all couplings of $X_{\lin m,\iy\lin}$ with a Markov chain $Z$
whose transition kernel is $M$. Then we have $\lim_{m\ri\iy}s_m=0$.
To be convinced of this convergence, consider for fixed $m\in\NN$, $(\wi X_n,\wi Y_n)_{n\in\NN}$ and $(\wi Z_n)_{n\in\NN}$
 two chains coupled as in the beginning of this section and starting from the initial conditions
$(\wi X_0,\wi Y_0)=(X_m,Y_m)$ and $\wi Z_0=X_m$.
Then (\ref{B}) and (\ref{ab2}) imply that the quantity $\sum_{n\in\NN}\EE[\vert \wi X_n-\wi Z_n\vert]$
converges exponentially fast to zero as $m$ goes to infinity.
\end{rem}

\appendix
\section{Numerical computations}

We present here the technical facts needed to conclude the proof of Proposition \ref{neigh}.
We used the  free numerical computational package Scilab 4.1.2. Next we investigate the Lipschitz constants alluded to in Remark \ref{LipDF}.
Finally we illustrate the results obtained in this paper by giving a numerical approximation
of 
the invariant measure $\mu$ of the chain $Z$ and of rate $L$ defined in (\ref{L}).\par\sm
But as a preliminary, let us give in Figure 1 the graphs of the six functions $z_i$, 
for $i\in\lin 1,6\rin$, since they are at the heart of this paper.

 \long\def\Checksifdef#1#2#3{
\expandafter\ifx\csname #1\endcsname\relax#2\else#3\fi}
\Checksifdef{Figdir}{\gdef\Figdir{}}{}
\def\dessin#1#2{
\begin{figure}[hbtp]
\begin{center}
\fbox{\begin{picture}(300.00,212.00)
\includegraphics{\Figdir}
\end{picture}}
\end{center}
\caption{\label{#2}#1}
\end{figure}}

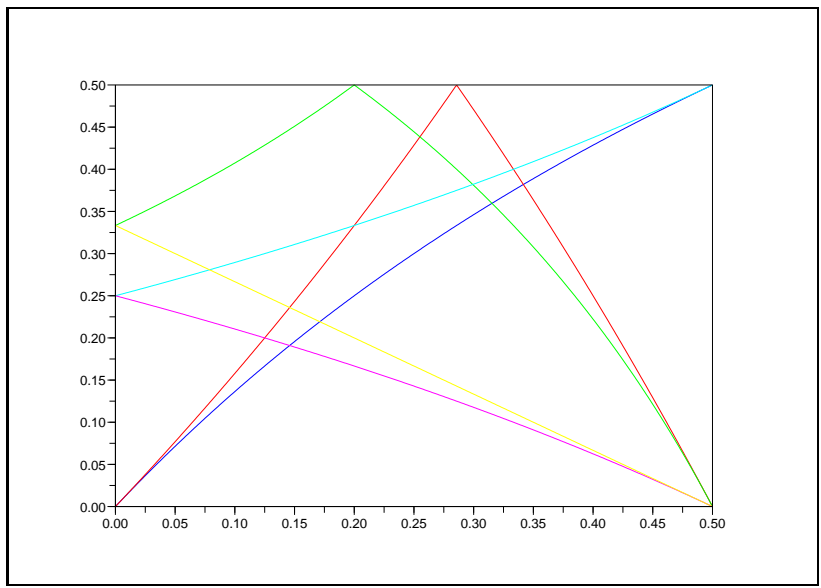
\dessin{Graphs of the $z_i$, for $i\in\lin 1,6\rin$}{fig1}

\subsection{Computations relative to Proposition \ref{neigh}}

The goal is to show that the function $F$ defined in  (\ref{F1}) does not vanish in $[0,1/2]$. 
Since we want to work formally as long as we can, it is better to consider
the piecewise rational fraction $R$ given in (\ref{R}).
Here are the codes we wrote
and the results given by the computer. \par\sm
$\bullet$ On the interval $(2/7,1/2)$.\par\sm 
First we compute the matrix $A\df(\sqrt{6}G(i,z_j(x))\sqrt{6}G(j,x))_{i,j\in\lin 1,6\rin}$,
where $x$ is a ``free variable". We multiplied the function $G$ by $\sqrt{6}$
to avoid square roots in the computations.
\begin{verbatim}
x=poly(0,'x');
z1=3/2*x/(1+x);
z2=(2-4*x)/(2-x);
z3=(2-4*x)/(3-3*x);
z4=(1+x)/(4-2*x);z5=(1-2*x)/(4-2*x);z6=(1-2*x)/3;
z=[z1 z2 z3 z4 z5 z6];
g1=2*(1+x);g2=(2-x);g3=3*(1-x);
g4=2*(2-x);g5=g4;g6=3;
g=[g1;g2;g3;g4;g5;g6];
a=horner(g,z);
A=a*diag(g)
\end{verbatim}
Running this code, we get the following polynonial matrix:
\bq
A&\df&\lt(\begin{array}{cccccc}
4 + 10x  &   8 - 10x   & 10 - 14x   &  10 - 2x   &   10 - 8x    & 8 - 4x\\
   4 + x    &   2 + 2x   &   4 - 2x   &   7 - 5x   &    7 - 2x   &   5 + 2x\\
    6 - 3x   &   9x      &    3 + 3x   &   9 - 9x   &    9        &   6 + 6x\\
     8 + 2x  &    4 + 4x  &    8 - 4x  &    14 - 10x &    14 - 4x  &  10 + 4x\\
     8 + 2x  &    4 + 4x  &    8 - 4x  &    14 - 10x &    14 - 4x  &  10 + 4x\\
     6 + 6x  &    6 - 3x  &    9 - 9x  &    12 - 6x  &    12 - 6x  &   9
\end{array}\rt)\eq
By construction, we have \bq
R(x)&=&\frac{1}{6^{36}}\prod_{i,j\in\lin 1,6\rin}A_{i,j}(x)\\
&=&\frac{2^{29}3^{17}}{6^{36}}\prod_{i,j\in\lin 1,6\rin}B_{i,j}(x)
\eq
where the matrix $B$ is defined by
\bq
B&\df&\lt(\begin{array}{cccccc}
  2 + 5x   &  4 - 5x  &  5 - 7x   &  5 - x    &  5 - 4x   &  2 - x\\
   4 + x   &    1 + x    &  2 - x    &  7 - 5x    &   7 - 2x    &  5 + 2x\\
     2 - x   &   x        &  1 + x   &   1 - x     & 1     &     1 + x\\
   4 + x   &  1 +x   &  2 - x   &   7 - 5x  &  7 - 2x   & 5 + 2x   \\
   4 + x   &   1 + x   &  2 - x   &   7 - 5x   &  7 - 2x   & 5 + 2x\\
    1 + x  &    2 - x  &    1 - x  &   2 - x   &   2 - x   &  1
\end{array}\rt)\eq
thus we obtain that $R(x)-1=P_1(x)$ on $[2/7,1/2]$, with
\bq
P_1(x)&\df&
2^{-7}3^{-19}(2 + 5x) (  4 - 5x )  (5 - 7x)   ( 5 - x )   (  5 - 4x )(  2 - x    )^8
   (4 + x ) ^3 ( 1 + x )^6\\&&(  7 - 5x  )^3(   7 - 2x  )^3 (  5 + 2x)^3x(  1 - x)^2-1\eq
Using the command \verb+roots+ of Scilab 4.1.2 which computes numerically the roots
of a given polynomial function, we obtained the roots of $P_1(x)$ depicted as crosses in Figure 2.
 
 \long\def\Checksifdef#1#2#3{
\expandafter\ifx\csname #1\endcsname\relax#2\else#3\fi}
\Checksifdef{Figdir}{\gdef\Figdir{}}{}
\def\dessin#1#2{
\begin{figure}[hbtp]
\begin{center}
\fbox{\begin{picture}(300.00,212.00)
\includegraphics{\Figdir}
\end{picture}}
\end{center}
\caption{\label{#2}#1}
\end{figure}}

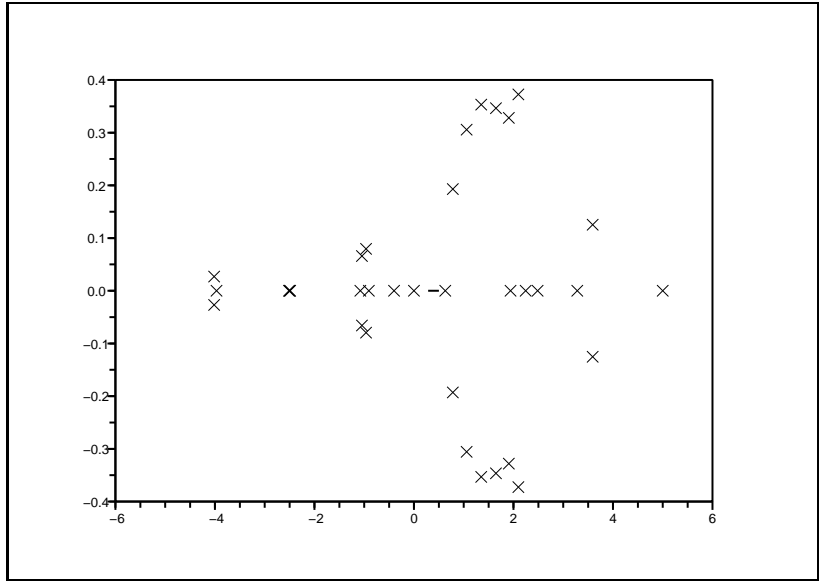
\dessin{Roots of $P_1(x)$}{fig2}
The segment $[2/7,1/2]$ is drawn as the
thick horizontal segment in the middle of the picture and it appears (by zooming) that 
none of the roots are on it. Indeed looking at the roots
found by the computer, the one on the left (respectively the right)
of the segment is 0.0006183 (resp.\ 0.6297289). The fact that $R>1$ on $[2/7,1/2]$
comes from the computation of $R(1/2)\approx 12.989284$. We also estimate that
$R(2/7)\approx 99.311045$. 
\par\me
$\bullet$ On the interval $(1/5,2/7)$.
\par\sm 
The first code above has only to be slightly modified: the definition of \verb+z2+
 (line 4)
is now
\begin{verbatim}
z2=(3*x)/(2-x);
\end{verbatim}
The corresponding matrix $A$ is also polynomial and the computer gives us
\bq
A&=&\lt(\begin{array}{cccccc}
   4 + 10x  &   4 + 4x  &   10 - 14x   &  10 - 2x    &  10 - 8x   &  8 - 4x    \\
    4 + x    &   4 - 5x   &   4 - 2x   &   7 - 5x    &   7 - 2x   &   5 + 2x    \\   
       6 - 3x   &   6 - 12x  &   3 + 3x   &   9 - 9x   &    9     &      6 + 6x   \\     
     8 + 2x  &    8 - 10x   &  8 - 4x  &    14 - 10x &    14 - 4x  &  10 + 4x    \\   
    8 + 2x   &   8 - 10x   &  8 - 4x  &    14 - 10x  &   14 - 4x  &  10 + 4x  \\  
        6 + 6x    &  6 - 3x   &   9 - 9x    &  12 - 6x   &   12 - 6x  &   9     
\end{array}\rt)\eq
Then we compute the roots of the polynomial function $P_2(x)$ representing $R(x)-1$
on $[1/5,2/7]$. As can be seen in Figure 3, none of them belong to this segment (drawn as thick line, the closest root is 0.4569663).
 
 \long\def\Checksifdef#1#2#3{
\expandafter\ifx\csname #1\endcsname\relax#2\else#3\fi}
\Checksifdef{Figdir}{\gdef\Figdir{}}{}
\def\dessin#1#2{
\begin{figure}[hbtp]
\begin{center}
\fbox{\begin{picture}(300.00,212.00)
\includegraphics{\Figdir}
\end{picture}}
\end{center}
\caption{\label{#2}#1}
\end{figure}}

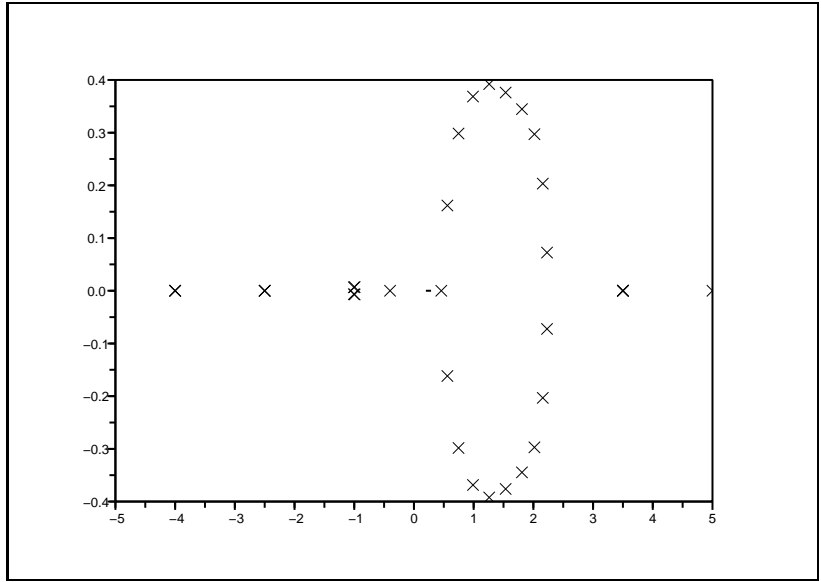
\dessin{Roots of $P_2(x)$}{fig3}
Thus $R$ is strictly positive on $[1/5,2/7]$, since we already knew that $R(2/7)>0$.
We estimate that $R(1/5)\approx 418.66239$.
\par\me
$\bullet$ On the interval $(0,1/5)$.
\par\sm 
The additional modification is the definition of \verb+z3+ (line 5),
\begin{verbatim}
z3=(1+x)/(3-3*x);
\end{verbatim}
Again we obtain a polynomial matrix $A$,
\bq
A&=&\lt(\begin{array}{cccccc}
   4 + 10x  &   4 + 4x   &   8 - 4x     & 10 - 2x    &  10 - 8x   &  8 - 4x    \\
     4 + x   &    4 - 5x   &   5 - 7x  &    7 - 5x   &    7 - 2x   &   5 + 2x    \\
   6 - 3x    &  6 - 12x  &   6 - 12x   &  9 - 9x     &  9       &    6 + 6x       \\
       8 + 2x   &   8 - 10x   & 10 - 14x  &   14 - 10x   &  14 - 4x   & 10 + 4x     \\
   8 + 2x   &   8 - 10x  &  10 - 14x   &  14 - 10x   &  14 - 4x  &  10 + 4x   \\
       6 + 6x    &  6 - 3x   &   9 - 9x    &  12 - 6x    &  12 - 6x   &  9  
\end{array}\rt)\eq
and Figure 4 gives a picture of the roots of 
$P_3(x)$ representing $R(x)-1$
on $[0,1/5]$. None of them belong to this segment (the closest one is 0.3995404).
 
 \long\def\Checksifdef#1#2#3{
\expandafter\ifx\csname #1\endcsname\relax#2\else#3\fi}
\Checksifdef{Figdir}{\gdef\Figdir{}}{}
\def\dessin#1#2{
\begin{figure}[hbtp]
\begin{center}
\fbox{\begin{picture}(300.00,212.00)
\includegraphics{\Figdir}
\end{picture}}
\end{center}
\caption{\label{#2}#1}
\end{figure}}

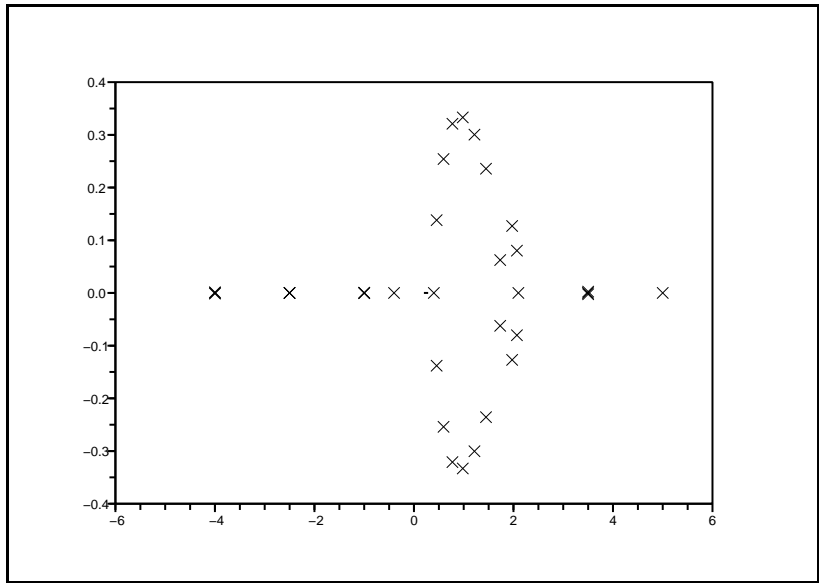
\dessin{Roots of $P_3(x)$}{fig4}
Thus $R$ is also strictly positive on $[0,1/5]$, since we have already seen that $R(1/5)>0$.
We estimate that $R(0)\approx 13496.561$.
\par\sm\sm
In conclusion, the function $F$ defined in (\ref{F1}) is positive on $[0,1/2]$,
 it is illustrated in Figure~5.
 
 \long\def\Checksifdef#1#2#3{
\expandafter\ifx\csname #1\endcsname\relax#2\else#3\fi}
\Checksifdef{Figdir}{\gdef\Figdir{}}{}
\def\dessin#1#2{
\begin{figure}[hbtp]
\begin{center}
\fbox{\begin{picture}(300.00,212.00)
\includegraphics{\Figdir}
\end{picture}}
\end{center}
\caption{\label{#2}#1}
\end{figure}}

\dessin{Graph of function $F$}{fig5}
There is a more analytical approach to justify the positivity of $F$. Indeed, the above polynomial expressions
for the three matrices $A$ (one on each of the segments $[0,1/5]$, $[1/5,2/7]$
 and $[2/7,1/2]$)
gives a simple uniform bound $C$ on $F'$ on $\in(0,1/2)\setminus\{1/5, 2/7\}$ via the formula
\bq
\fo x\in(0,1/2)\setminus\{1/5, 2/7\},\qquad
 \lve F'(x)\rve&=& \frac1{36}\lve \sum_{i,j\in\lin 1,6\rin}\frac{A_{i,j}'(x)}{A_{i,j}(x)}\rve\\
 &\leq&\max_{i,j\in\lin 1,6\rin}\lve \frac{A_{i,j}'(x)}{A_{i,j}(x)}\rve
 \eq
 Working out the details, we get that we can take $C=14/3$ (obtained for $i=3$, $j=2$ and
 $x=2/7$ for the matrix $A$ valid on the interval $[1/5,2/7]$).
Thus if there exists $N\in\NN^*$ such that 
\bq
\min_{i\in\lin 0,N\rin}F(i/N)-14/(6N)&>&0\eq
we can conclude to the positivity of $F$ on $[0,1/2]$. We tried with 
$N=\lceil 14/(3F(1/2))\rceil=66 $ and it works (except that we have not taken into account the precision of the computer).\par\me
The fact that $F$ and $R$ are decreasing on [0,1/2] can also be checked numerically.
Indeed, using the command \verb+derivat+ of Scilab 4.1.2, we
can compute first the derivatives of the polynomial functions $P_1$, $P_2$
and $P_3$ and next their roots. It appears that they do not belong
respectively to the segments $[2/7,1/2]$, $[1/5,2/7]$
 and $[0,1/5]$. Since we have already seen that $F(0)>F(1/5)>F(2/7)>F(1/2)$, it follows that
$F$ is decreasing on $[0,1/2]$ and thus bounded below by $F(1/2)$.\par\sm
As suggested to us by
André Galligo,
there is also a totally algebraic approach to this problem via the use of Sturm's sequences.
Since we are dealing with polynomial functions with rational coefficients, this method
can be made rigorous, if we are allowed to work on the computer with numbers with a large
number of digits (here of order $42=\lceil\log_{10}(36 !)\rceil$, since we have to consider the values of
all the derivatives of the polynomial function of interest at the boundary of the interval
where we are looking for roots), in order to deal only
with true ``integers''. So this approach would require working with  mathematical software such as Maple.
\par
\me
Finally let us also draw the graph of the function 
\bq
\wi F\st[0,1/2]\ni x&\mapsto&\frac16\sum_{i\in\lin 1, 6\rin}\ln(G(i,x))\eq
to show that it was necessary to iterate the kernel $M$ 
to get a lower bound as the one in Proposition~\ref{neigh} (or
to be able to apply the result of Barnsley and Elton \cite{MR932532}). 
 
 \long\def\Checksifdef#1#2#3{
\expandafter\ifx\csname #1\endcsname\relax#2\else#3\fi}
\Checksifdef{Figdir}{\gdef\Figdir{}}{}
\def\dessin#1#2{
\begin{figure}[hbtp]
\begin{center}
\fbox{\begin{picture}(300.00,212.00)
\includegraphics{\Figdir}
\end{picture}}
\end{center}
\caption{\label{#2}#1}
\end{figure}}

\dessin{Graph of function $\wi F$}{fig6}
Indeed, $\wi F$ is negative in a neighborhood of $1/2$,
which means that the isoperimetric value of the triangles near this neighborhood
have a tendency to decrease. But fortunately, the abscissa of the characteristic point of new triangle obtained by the barycentric subdivision has a chance of 4
over 6 (see Figure 1)  to be close to zero, a
zone where the isoperimetric functional has a strong tendency to increase.

\subsection{On the Lipschitz constants of Remark \ref{LipDF}}

We implement on the computer the verification of the criterion (\ref{DF}) 
from Diaconis and Freedman \cite{MR1669737} in our particular situation.
We adopt an analytical approach, despite the fact we are essentially dealing
with polynomial functions.
Let us denote (recall the notation of (\ref{G}))
\bq
 \fo i\in\lin 1,6\rin,\,\fo x\in[0,1/2],\qquad
 D(i,x)&\df&\lve z'_i\rve(x)\ =\ \frac1{G^2(i,x)}\eq
 and 
 for any $n\in\NN^*$ and $(i,j)\in\lin 1,6\rin^2$, 
 \bq
 K_{i,j}(n)&\df&\sup_{x\in S_n} D(i,z_j(x))D(j,x)\eq
 where $S_n\df\lin 0, n\rin/(2n)$ is a mesh of $[0,1/2]$.
 Then for any $(i,j)\in\lin 1,6\rin^2$, the Lipschitz constant of $z_i\circ z_j$
 is given by
 \bq
 K_{i,j}&\df&\sup_{x\in [0,1/2]} D(i,z_j(x))D(j,x)\\
 &=&\sup_{n\in\NN^*} K_{i,j}(n)\eq
 Thus to show that criterion (\ref{DF}) is not satisfied by $M^2$,
 it is sufficient to check that for some $n\in\NN^*$, $\cK_2(n)\df\prod_{(i,j)\in\lin 1,6\rin^2}K_{i,j}(n)>1$.\\
Similarly for $M^3$, we define  
for any $n\in\NN^*$ and any $(i,j,k)\in\lin 1,6\rin^3$,
\bq K_{i,j,k}(n)
 &\df&\sup_{x\in S_n} D(i,z_j\circ z_k(x))D(j,z_k(x))D(k,x)\eq
so that
for any $(i,j,k)\in\lin 1,6\rin^3$, the Lipschitz constant of $z_i\circ z_j\circ z_k$
 is given by $K_{i,j,k}=\sup_{n\in\NN}K_{i,j,k}(n)$. 
We compute that
\bq
\max_{i\in\lin 1,6\rin,x\in[0,1/2]}D(i,x)&=&\frac83\\
\max_{i\in\lin 1,6\rin,x\in[0,1/2]}\lve\frac{dD(i,x)}{dx}\rve&=&\frac{32}3\eq 
so we easily deduce that
\bq
\max_{(i,j,k)\in\lin 1,6\rin^3,x\in[0,1/2]}\lve \frac{d}{dx}D(i,z_j\circ z_k(x))D(j,z_k(x))D(k,x)\rve
&\leq& \frac{195552}{243}\ \leq\ 805\eq
It follows that 
\bq
\fo n\in\NN^*,\,\fo (i,j,k)\in\lin 1,6\rin^3,\qquad K_{i,j,k}&\leq& K_{i,j,k}(n)+\frac{202}{n}\eq
Thus to see that criterion (\ref{DF}) is satisfied by $M^3$,
 it is sufficient to check that for some $n\in\NN^*$, \bq\cK_3(n)&\df&\prod_{(i,j,k)\in\lin 1,6\rin^3}\lt(K_{i,j,k}(n)+\frac{202}{n}\rt)<1\eq
The following program (where the index set $\lin 1,6\rin$ has been replaced by the more convenient set
$\lin 0,5\rin$) gives that
 $\cK_2(1)=\verb+k2+\approx 92.067225$ and $\cK_3(10000)=\verb+k3+\approx0.0116774$, values which justify the assertions made
 in Remark \ref{LipDF}.
 \begin{verbatim}
function z=Z(x,i)
select i
case 0 then z=3/2*x/(1+x);
case 1 then z=3*bool2s((x<2/7))*x/(2-x)+bool2s((x>=2/7))*(2-4*x)/(2-x);
case 2 then z=bool2s((x<1/5)).*(1+x)./(3-3*x)+bool2s((x>=1/5)).*(2-4*x)./(3-3*x);
case 3 then z=(1+x)/(4-2*x);
case 4 then z=(1-2*x)/(4-2*x);
case 5 then z=(1-2*x)/3;
end; endfunction

function d=D(x,i)
select i
case 0 then d=3/2 ./(1+x).^2;
case 1 then d=6 ./(2-x).^2;
case 2 then d=2/3 ./(1-x).^2;
case 3 then d=3/2 ./(2-x).^2;
case 4 then d=3/2 ./(2-x).^2;
case 5 then d=2/3
end; endfunction

function d=DD(x,i)
j=floor(i/6); k=i-6*j;
d=D(Z(x,k),j)*D(x,j);
endfunction

function d=DDD(x,i)
j=floor(i/36); k=floor((i-36*j)/6); l=i-36*j-6*k ;
d=D(Z(Z(x,l),k),j)*D(Z(x,l),k)*D(x,l)+202/10000;
endfunction

r2=[0 1/2]; I2=0:35; K2=feval(r2,I2,DDD);
k2=prod(max(K2,'r'))

r3=0:1/(2*10000):.5; I3=0:215; K3=feval(r3,I3,DDD);
k3=prod(max(K3,'r'))
\end{verbatim}

\subsection{Approximation of the invariant measure $\mu$}

Proposition \ref{slln} leads to an approximation of the invariant measure $\mu$ through the strong law of large numbers. 
 In practice, we run the following code, which simulates the chain $(Z_n)_{n\in\lin 1,N\rin}$
 with $N=100000$ starting from an initial variable $Z_1$ chosen uniformly over $[0,1/2]$.
  \begin{verbatim} 
n=100000;
z(1)=.5*rand(); s(1)=0;
for i=2:n
j=floor(6*rand());
z(i)=Z(z(i-1),j);
s(i)=s(i-1)-log(D(z(i-1),j))/2;
end
L=s(n)/n 
  \end{verbatim}
 The program takes into account the functions \verb+Z+  and \verb+D+      defined
 in the code of the previous subsection.
 Figure 7 shows the approximation we got of $\mu$ through a histogram using 100 bars.
 
 \long\def\Checksifdef#1#2#3{
\expandafter\ifx\csname #1\endcsname\relax#2\else#3\fi}
\Checksifdef{Figdir}{\gdef\Figdir{}}{}
\def\dessin#1#2{
\begin{figure}[hbtp]
\begin{center}
\fbox{\begin{picture}(300.00,212.00)
\includegraphics{\Figdir}
\end{picture}}
\end{center}
\caption{\label{#2}#1}
\end{figure}}

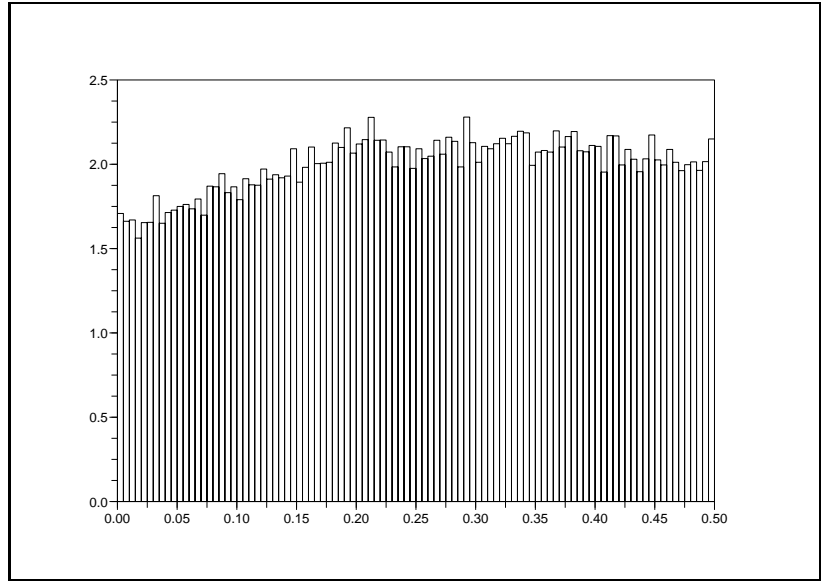
\dessin{An approximation of $\mu$}{fig7}
Finally, the code provides an approximation of the rate $L$ defined in (\ref{L}).
 Running the program several times, we always got values between 0.07 and 0.08, suggesting that
 Theorem \ref{th1} should be satisfied with $\chi=0.07$.

\bigskip
\par\hskip5mm\textbf{\large Acknowledgments:}\par\sm\noindent 
We are indebted to David Blackwell, Curtis~T. McMullen and Bob Hough for kindly providing us
their manuscripts \cite{Blackwell_barycentric,Diaconis_barycentric,MR2516262} in preparation on the subject.
We acknowledge the generous support of the ANR's Chaire  d'Excellence
program of the CNRS. We are very grateful to Oliver Riordan and to the referees which helped us to 
improve the paper by simplyfing its arguments. In particular the suggestion to use the simpler functional $J$
instead of $I$ in Section 2 is due to one of them.

\def\cprime{$'$}

\vskip15mm
\hskip10mm\box6

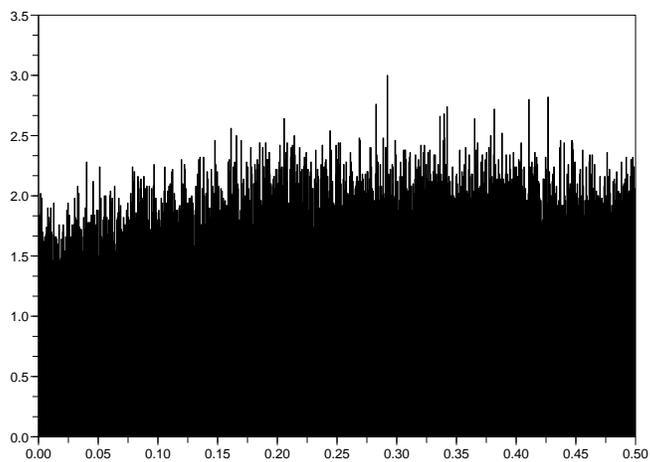
\end{document}